\documentclass[11pt,a4paper]{article}
\usepackage[margin=1in]{geometry}
\usepackage[T1]{fontenc}
\usepackage{lmodern}
\usepackage{amsmath,amssymb,amsthm}
\usepackage{enumitem,microtype}
\usepackage[hidelinks,unicode]{hyperref}
\usepackage[capitalise,nameinlink,noabbrev]{cleveref}

\theoremstyle{plain}
\newtheorem{theorem}{Theorem}[section]
\newtheorem{lemma}[theorem]{Lemma}
\newtheorem{proposition}[theorem]{Proposition}
\newtheorem{corollary}[theorem]{Corollary}

\newcommand{\br}{\operatorname{br}_3}
\newcommand{\eps}{\varepsilon}
\newcommand{\ind}[1]{\mathbf{1}_{#1}}
\newcommand{\dotcup}{\mathbin{\dot\cup}}
\newcommand{\Z}{\mathbb{Z}}
\setlist[enumerate]{label=\textup{(\roman*)},leftmargin=*,itemsep=3pt}

\title{The three-colour bipartite Ramsey number 
of sufficiently long even cycles}

\author{Jiangdong Ai\thanks{School of Mathematical Sciences and LPMC, Nankai University,
Tianjin 300071, China. Email: \texttt{jd@nankai.edu.cn}. Partially supported by the
National Natural Science Foundation of China (No. 12522117), the Fundamental and Interdisciplinary
Disciplines Breakthrough Plan of the Ministry of Education of China (JYB2025XDXM207) and
the Fundamental Research Funds for the Central Universities, Nankai University.}
\and
Yihan Yin\thanks{School of Mathematical Sciences and LPMC, Nankai University,
Tianjin 300071, China. Email: \texttt{2549152459@qq.com}.}}

\date{}

\hypersetup{
  pdftitle={The three-colour bipartite Ramsey number of sufficiently long even cycles},
  pdfauthor={Jiangdong Ai and Yihan Yin},
  pdfsubject={Bipartite Ramsey theory},
  pdfkeywords={bipartite Ramsey number, even cycle, connected matching, vertex cover}
}

\begin{document}
\maketitle

\begin{abstract}
For a bipartite graph $H$, let $\operatorname{br}_3(H)$ be the least
integer $N$ such that every three-colouring of $K_{N,N}$ contains
a monochromatic copy of $H$. We prove that
$\operatorname{br}_3(C_{2n})=3n$ for all sufficiently large $n$,
and give the lower bound for every $n\ge2$. For the upper bound,
we reduce a counterexample to a finite rational weighted profile.
Cycle extension and vertex-cover arguments force the critical
components to have six unique minimum covers which partition the
profile. The remaining component supports give a contradiction.


\end{abstract}

\section{Introduction}\label{sec:1}

For bipartite graphs $H_1,\ldots,H_r$, let
$\operatorname{br}(H_1,\ldots,H_r)$ be the least positive integer
$N$ such that every $r$-colouring of $E(K_{N,N})$ contains a copy
of $H_i$ in colour $i$ for some $i\in[r]$. When all graphs equal
$H$, write $\operatorname{br}_r(H)$; for three colours we use
$\br(H)$. As usual, $P_k$ and $C_k$ denote a path and a cycle on
$k$ vertices.

Gy\'arf\'as and Lehel~\cite{GL1973} and, independently, Faudree
and Schelp~\cite{FS1975} initiated the study of bipartite Ramsey
numbers. Their path results give
\[
 \operatorname{br}(P_{2m},P_{2n})=m+n-1\qquad(m,n\ge2).
\]
For even cycles, Zhang and Sun~\cite{ZS2011}, Zhang, Sun and
Wu~\cite{ZSW2013}, and Gholami and Rowshan~\cite{GR2023}
determined the two-colour numbers when one cycle is $C_4$, $C_6$
or $C_8$, respectively. The small-bandwidth result of Shen, Lin
and Liu~\cite{SLL2018} implies
$\operatorname{br}_2(C_{2n})=(2+o(1))n$.
Yan and Peng~\cite{YP2024} settled the remaining cases: for
$n\ge m\ge5$,
\[
 \operatorname{br}(C_{2n},C_{2m})=
 \begin{cases}
 n+m-1,&n>m,\\
 2n,&n=m.
 \end{cases}
\]
Thus the diagonal cycle number is one larger than the
corresponding path number.

Buci\'c, Letzter and Sudakov~\cite{BLS2019} proved that
\[
 \br(C_{2n})=(3+o(1))n,
 \qquad \br(P_{2n})=(3+o(1))n.
\]
They first determined the relevant connected-matching Ramsey
numbers and then applied regularity. Luo and Peng~\cite{LuoPeng2020}
extended this approach to the asymmetric three-colour problem.
For fixed $r\ge3$, Liu and Peng~\cite{LiuPeng2024} showed that,
if positive reals $\alpha_1,\ldots,\alpha_r$ satisfy
\[
 \alpha_{j+2}\ge\bigl((j+2)!-1\bigr)
 \sum_{i=1}^{j+1}\alpha_i\qquad(1\le j\le r-2),
\]
then
\[
 \operatorname{br}(C_{2\lfloor\alpha_1n\rfloor},\ldots,
 C_{2\lfloor\alpha_rn\rfloor})
 =\left(\sum_{i=1}^r\alpha_i+o(1)\right)n.
\]
The diagonal problem behaves differently for more colours.
Constructions of DeBiasio, Gy\'arf\'as, Krueger, Ruszink\'o and
S\'ark\"ozy~\cite{DGKRS2020} give
$\operatorname{br}_4(C_{2n})\ge(5-o(1))n$ and
$\operatorname{br}_r(C_{2n})\ge(2r-4-o(1))n$ for $r\ge5$.
Buci\'c, Letzter and Sudakov~\cite{BLSmulti2019} proved
$\operatorname{br}_4(C_{2n})=(5+o(1))n$ and the upper bound
$(2r-3+o(1))n$ for fixed $r\ge4$, with further improvements for
larger $r$ and analogous path results.

A \emph{connected matching} is a matching contained in one
component. \L uczak's method~\cite{Luczak1999} reduces asymptotic
cycle problems to connected matchings in almost complete reduced
graphs. Letzter~\cite{Letzter2022} gave a further reduction to
complete blow-ups, including bipartite hosts. These methods give
the leading term, but an exact cycle number also requires control
of the vertices lost in the reduction.

For three colours, the construction in~\cite{BLS2019} gives
$\br(C_{2n})\ge3n-2$, while the upper bound there is $3n+o(n)$.
We determine the exact value for sufficiently large $n$.

\begin{theorem}\label{r1.1}
There exists an integer $n_0$ such that
\[
 \br(C_{2n})=3n
 \qquad\text{for every integer }n\ge n_0.
\]
\end{theorem}

We first give a colouring which proves the lower bound for every
$n\ge2$. We then describe the reduction used for the upper bound.

\subsection{A lower bound}

Put $m=n-1$. Partition the two shores into
\[
 L=A_0\dotcup A_1\dotcup A_2\dotcup\{x,z\},
 \qquad
 R=B_0\dotcup B_1\dotcup B_2\dotcup\{p,q\},
\]
where $|A_i|=|B_i|=m$ for $i\in\{0,1,2\}$. Colour every edge
between a row set and a column set according to the following
table. We denote this colouring by $F_n$.
\[
\begin{array}{c|ccccc}
 &B_0&B_1&B_2&p&q\\\hline
A_0&0&1&2&1&2\\
A_1&1&2&0&1&1\\
A_2&2&0&1&2&2\\
x&0&0&0&0&0\\
z&0&2&1&1&2
\end{array}
\]

\begin{proposition}\label{r1.2}
For every $n\ge2$, the colouring $F_n$ has no monochromatic
$C_{2n}$. Consequently, $\br(C_{2n})\ge3n$.
\end{proposition}
\begin{proof}
The colour-$0$ edges form
\[
 K_{A_0\cup\{x,z\},B_0},\quad
 K_{A_1\cup\{x\},B_2},\quad
 K_{A_2\cup\{x\},B_1},
\]
which meet only at $x$, together with the bridges $xp,xq$.
The colour-$1$ edges form
\[
 K_{A_0,B_1\cup\{p\}},\quad
 K_{A_1,B_0\cup\{p,q\}},\quad
 K_{A_2\cup\{z\},B_2},
\]
together with the bridge $zp$. The first two graphs meet only
at $p$; the third is attached by $zp$. Every cycle in either
colour therefore lies in one displayed graph, whose smaller
shore has size $m$.
Interchanging $A_1,A_2$, $B_1,B_2$, $p,q$, and colours $1,2$
preserves the table, so every monochromatic cycle has length at
most $2m$. Both shores have $3m+2=3n-1$ vertices, proving the
lower bound.
\end{proof}

\subsection{An overview of the upper bound}

From a counterexample on $K_{3n,3n}$, regularity and rational
correction give a finite weighted profile. Its shores have
weight $3h=1$, and every monochromatic component has a vertex
cover of weight at most $h$; equality defines a \emph{critical}
component. The finite alternative in \cref{struct:alternative}
gives a wide common strip, a three-row permutation profile, or
a component with shore weights $h$ and at least $2h$. Exact
cycle constructions exclude all three in the actual graph.

The main tools are the block--cutpoint structure of Latin cores
and an integer cover-trade identity. The latter forces two
critical components in each colour, whose six unique minimum
covers partition the profile. Rectangles between these covers
then contradict the remaining component supports.

In \cref{sec:2} we develop extraction and cover trades;
\cref{sec:3,sec:permutations} give the Latin-core, permutation
and transfer arguments. We prove the width and tripartition
obstructions in \cref{sec:obstructions} and the critical-cover
partition in \cref{sec:critical}. The residual contradiction and
the main proof are in \cref{sec:residual}.

We also call colours $0,1,2$ red, blue and green. The integer $n$
is the target half-length of an actual cycle. Profile shores
have weight one until \cref{sec:13}, where weights are scaled
by three.

\section{Finite profiles and exact extraction}\label{sec:2}

A \emph{profile} has finite shores $X,Y$, positive rational
vertex weights summing to one on each shore, and nonempty
allowed-colour sets $\mathcal A(x,y)\subseteq\{0,1,2\}$.
Put $h=1/3$, $w(S)=\sum_{v\in S}w_v$, and
$W_c(x,y)=\ind{c\in\mathcal A(x,y)}/|\mathcal A(x,y)|$.
Components and covers refer to the allowed-edge graphs.
A component is \emph{critical} if its minimum cover weight is
$h$, and \emph{pure-sided} if each of its minimum covers is
contained in a single shore. Sets consist of whole vertices; degrees are finite sums,
such as $d_c(y,S)=\sum_{x\in S}w_xW_c(x,y)$.
We transfer only absence, presence and entirely one-colour
rectangles; the values $W_c$ need not approximate actual densities.

Our inputs from~\cite{BLS2019} are Proposition~6,
$r(k,k+1,k+1)\le3k$, simultaneous regularity (Lemma~14), and
connected-matching lifting (Lemma~16). We apply Lemma~16 to each
colour class separately, so the colour-$c$ reduced graph has an edge
for every pair that is $\eps$-regular in colour $c$ with colour-$c$
density at least $d$. A matching of size $s$ in one component of
this graph, with $d>20\eps$ and cluster size $\ell$, gives all even
colour-$c$ cycle lengths up to $2(1-9\eps/d)s\ell$.

\begin{lemma}\label{finite:correction}
Consider a finite family of rational linear systems on a fixed
bounded rational box. There are constants $a,L>0$ such that every
rational point in the box whose maximum constraint violation is
$\eta<a$ is within maximum-norm distance $L\eta$ of a rational
solution of the corresponding system.
\end{lemma}
\begin{proof}
Let $R(x)$ be the maximum of the absolute equation residuals
and positive inequality residuals. Minimising $R$ on the box is
a rational linear programme; bound the auxiliary residual
variable above all residuals on the box. Each infeasible system
has positive rational optimum, so choose $a$ below all such
optima in the finite family.

For a feasible system, let $P$ be its solution polytope in the
box and $D(x)=\min_{y\in P}\|x-y\|_\infty$. Enumerating
linear-programme bases gives a finite rational subdivision on
which both $R$ and $D$ are piecewise affine. On a common
refinement choose $L$ at least $D/R$ at all vertices with
$R>0$; where $R=0$, also $D=0$. Interpolation gives $D\le LR$.
Rational inputs admit rational optimal solutions. Take common
constants over the finite family.
\end{proof}

\begin{theorem}\label{r2.2}
Put $K=11^3$. For every sufficiently small positive function $F$,
there are scales $t_0>\cdots>t_{M+1}>0$, where $M=M(K)$ and
$t_{i+1}\le F(t_i)$, with the following property. For all
sufficiently large $n$, every three-colouring of $K_{3n,3n}$ with
no monochromatic $C_{2n}$ gives a profile on at most $K$ vertices
per shore and consecutive scales $t=t_j,s=t_{j+1}$ such that:
\begin{enumerate}
\item every monochromatic component has a cover of weight at most $h$;
\item for every $A\subseteq X\dotcup Y$ and
$r\in\{0,h,2h,1,2\}$, either $w(A)=r$ or $|w(A)-r|\ge t/2$;
\item there are disjoint actual vertex sets $V_v$ satisfying
$\bigl||V_v|/(3n)-w_v\bigr|\le C_Ks$, with at most $C_Ksn$
vertices omitted from each shore;
\item on $V_x\times V_y$, colour $c$ has at least $tn^2$ edges if
it is allowed, and at most $sn^2$ edges otherwise.
\end{enumerate}
In particular, every positive vertex weight is at least $t/2$.
The scales are fixed before the regularity parameters and $n$.
Here and below, $C_K$ denotes a sufficiently large constant
depending only on $K$.
\end{theorem}
\begin{proof}
On $[0,1]^{2K}$ pre-enumerate shore totals, nonnegativity,
every possible choice of at most fifteen marked-cover
inequalities, and every selection of equations $w(A)=r$,
where $A\subseteq X\dotcup Y$ and $r\in\{0,h,2h,1,2\}$.
Let $a,L$ be the correction constants. Choose $t_0<a$ and
the hierarchy so that $2KLF(t)<t/2$ and $3F(t)<9t^2$.

Choose $\eta\ll t_{M+1}$ and then $d,\eps/d\ll\eta$.
Apply simultaneous regularity, retaining each regular pair
of density at least $d$ in each colour. By lifting, a retained
component has no matching of size $n/((1-9\eps/d)\ell)$.
K\H{o}nig's theorem gives a cluster cover of normalised weight
at most $h+\eta$. The total discarded edge count, including
irregular pairs, low-density colours and exceptional vertices,
is at most $\eta n^2$.

Merge colour components while their combined two-shore mass is
at most $2h$. A merged block is covered by the smaller shores
of its original components, of total weight at most $h$;
an unmerged block keeps its cover. At termination, any two
blocks have combined mass greater than $2h$. Six blocks,
paired into three pairs, would exceed the total mass two.
Thus each colour has at most five blocks. Refining by block
and cover membership gives at most $(1+5\cdot2)^3=K$ types
per shore. Include exceptional vertices in outside types
and pad with zero coordinates.

Write $u_v=|V_v|/(3n)$ and $z_{cxy}=e_c(V_x,V_y)/n^2$.
Marks permit colour $c$ only on pairs inside one $c$-block (not
the outside label) that meet its cover. Unpermitted edges have total
normalised mass at most $\eta$. There are at most $M(K)$
numbers $z_{cxy}$ and $|u(A)-r|$; we may take
$M=3K^2+5\cdot2^{2K}$. One interval $(s,t)$ in the hierarchy
contains none. Impose all equations with residual at most
$s$, shore totals, nonnegativity and marked-cover bounds $h$.
The vector $u$ violates this system by at most $s$, so
correction gives rational $w$ with $\|w-u\|_\infty\le Ls$.

Each unselected affine form $u(A)-r$ changes by less than
$t/2$, hence retains its sign and has magnitude at least
$t/2$. Singleton equations make $w_v=0$ when $u_v\le s$;
all other types have $u_v\ge t$ and remain positive. Delete
the zero types, omitting actual normalised mass at most $2Ks$.

Set $\mathcal A(x,y)=\{c:z_{cxy}\ge t\}$ on surviving pairs.
It is nonempty, since otherwise
$9u_xu_y=\sum_cz_{cxy}\le3s<9t^2$, contrary to $u_x,u_y\ge t$.
As $\eta<s$, every allowed colour is permitted by its marks.
Its components lie in individual marked blocks and inherit
their corrected covers. The correction and the empty interval
give the remaining bounds.
\end{proof}

\begin{proposition}\label{r2.3}
If every monochromatic component has a cover of weight at most
$h$, then each colour has at least one and at most three critical
components.
\end{proposition}
\begin{proof}
If one colour has no critical component, choose strict covers
in that colour and covers of weight at most $h$ in the other
two. Clear denominators, replacing each vertex $v$ by $3kw_v$
copies, so that each shore has size $3k$. Choose any allowed
colour on each copy pair. Components in the first colour
have covers of size at most $k-1$; in the others they have
covers of size at most $k$, contradicting
$r(k,k+1,k+1)\le3k$.
Each critical component has weight at least $h$ on both
shores, so at most three can be disjoint in one colour.
\end{proof}

\begin{lemma}\label{finite:trade}
Let $H_1,\ldots,H_m\subseteq X\dotcup Y$, where $X,Y$ are finite
disjoint shores, and let $\varnothing\ne I\subseteq[m]$. Exactly
one of the following holds:
\begin{enumerate}
\item there is an integer-valued vector $z$ such that
$z(X)=z(Y)=0$, $z(H_i)\le0$ for every $i$, and $z(H_i)<0$ for
every $i\in I$;
\item there are nonnegative integers $a_i,k_X,k_Y$ such that
\[
 \sum_{i\in I}a_i>0,\qquad
 \sum_i a_i\ind{H_i}=k_X\ind X+k_Y\ind Y.
\]
\end{enumerate}
\end{lemma}
\begin{proof}
Scalar products show that the alternatives are incompatible.
Eliminate the shore equalities, then combine each upper bound
with each lower bound using positive integer multipliers, marking a row
strict if any source is strict. If no contradiction occurs,
rational back-substitution and clearing denominators give the
trade. Otherwise the combination yielding $0<0$ uses a strict
row. Restoring the shore equalities gives the incidence
identity, whose shore constants are nonnegative because its
left-hand side is coordinatewise nonnegative.
\end{proof}

\begin{lemma}\label{r2.4}
Suppose that every monochromatic component has a cover of weight
at most $h$. Choose an arbitrary minimum cover $H_{cj}$ for every
critical component. There are nonnegative integers $a_{cj},k_X,k_Y$ such
that
\[
 \sum_ja_{cj}>0\quad(c=0,1,2),\qquad
 \sum_{c,j}a_{cj}\ind{H_{cj}}
 =k_X\ind X+k_Y\ind Y.
\]
\end{lemma}
\begin{proof}
Fix a colour $c$. A sufficiently small rational multiple of
a trade strictly decreasing all its chosen critical covers
would preserve positive weights, shore totals and the other
critical-cover bounds. Fixed covers of noncritical components
would remain strict. Colour $c$ would then have no critical
component, contrary to \cref{r2.3}. Apply \cref{finite:trade}
with the colour-$c$ covers as the strict index set, and sum
the identities for the three colours. This auxiliary change
of weights is used only in the finite profile.
\end{proof}

\begin{lemma}\label{r2.12}
For each colour $i$, let a colour-$i$ component have whole left
shore $A_i$ and whole right shore $B_i$, where $A_0,A_1,A_2$
partition $X$ and $w(A_i)=h$. Put $D=\bigcap_iB_i$. Then every pairwise
intersection of the $B_i$ equals $D$, every $A_i\times D$ is
entirely colour $i$, and
\[
 w(D)\ge\frac{\sum_iw(B_i)-3h}{2}.
\]
\end{lemma}
\begin{proof}
For distinct $i,j,k$, colours $i,j$ are forbidden on
$A_k\times(B_i\cap B_j)$. Thus $B_i\cap B_j\subseteq B_k$
and this rectangle has only colour $k$. All pairwise
intersections equal $D$, and $A_i\times D$ has colour $i$.
The bound follows by inclusion--exclusion and $w(Y)=3h$.
\end{proof}

The following finite result will complete the reduction. We
state it here and prove it in \cref{sec:structural-proof}.

\begin{theorem}\label{struct:alternative}
Let a finite rational profile have shores of weight $3h$, and
suppose every monochromatic component has a vertex cover of
weight at most $h$. After interchanging the shores if necessary,
at least one of the following holds:
\begin{enumerate}
\item there is a partition $X=F_0\dotcup F_1\dotcup F_2$ with
$w(F_i)=h$, and a set $D\subseteq Y$ with $w(D)>h$, such that
$F_i\times D$ is entirely colour $i$ for each $i$;
\item there is a partition $X=F_0\dotcup F_1\dotcup F_2$ with
$w(F_i)=h$ such that every $y\in Y$ has a permutation
$\pi_y\in S_3$ for which $F_i\times\{y\}$ is entirely colour
$\pi_y(i)$;
\item some monochromatic component has whole-shore weights
$h$ and at least $2h$.
\end{enumerate}
\end{theorem}

\section{Homogeneous Latin cores}\label{sec:3}

Let $A_i,B_j$ be large classes of sizes $n-a_i,n-b_j$, with
colour $i+j$ on $A_i\times B_j$. Indices and colours are in
$\Z/3\Z$. Each exceptional vertex has a constant colour towards
each opposite large class; the corner $L_*\times R_*$ is arbitrary.
Write
\[
 r_L=|L_*|=\sum_i a_i,\quad r_R=|R_*|=\sum_jb_j,
 \quad r=r_L+r_R,\qquad n\ge3r+6.
\]
A \emph{target} is a monochromatic $C_{2n}$. In a counterexample,
all deficits are positive. Indeed, a class of size at least $n$
has three common colour neighbourhoods covering the opposite
$3n$ vertices, one of which gives a $K_{n,n}$.

\begin{lemma}\label{r3.1}
An exceptional vertex with a large-class vector may be recoloured
at the corner and added to that class. Every target obtained
afterwards can be recovered before the operation.
\end{lemma}
\begin{proof}
For a left absorption into a class of size $m$, put $q=r_R$.
Then $m+1>2q$. On a new target, replace the added vertex by an
unused old class vertex if possible. Otherwise all $m+1$ vertices
occur, and at most $2q$ of their positions touch right exceptional
vertices. Permute the class vertices so that the added vertex
avoids every recoloured edge. Each operation decreases $r$ and
one deficit. If a deficit becomes zero, recover a target from
the resulting $K_{n,n}$ by undoing the operations.
\end{proof}

After these operations, the model is \emph{normalised}. Its
remaining permutation vectors are $X_s(j)=s-j$ on the left and
$Y_t(i)=t-i$ on the right. Other vectors are \emph{ordinary},
with multiplicities $(2,1,0)$, or constant; both are called
\emph{repeated}. The principal colour-$c$ core of $X_s$ is
\[
 Q_c(s)=(i,j)=(-c-s,s-c).
\]
The types $X_s,Y_t$ are compatible exactly when $s+t=0$, in which
case their principal cores agree in every colour. The structural
arguments below concern a normalised counterexample.

\begin{lemma}\label{r3.2}
Let $p_s=|X_s|$. For any integer
\[
 0\le h\le\min\{\min_s p_s,\min_i a_i\},
\]
we may subtract $h$ from every $p_s,a_i$, preserving recovery,
core separators and activity (defined below). Consequently, $\min_i a_i>\min_s p_s$, and one
permutation type may be made empty. The assertions transpose.
\end{lemma}
\begin{proof}
Choose $h$ vertices of each $X_s$, add $h$ chosen vertices to
each $A_i$, and recolour accordingly. On a new colour-$c$ target,
the excess $t_i\le h$ over the old $|A_i|$ may occupy positions
with two large-class neighbours, since $|A_i|>2r_R$. These
neighbours lie in $B_{c-i}$, so fill the positions with chosen
$X_{-c-i}$ vertices. Different rows use different chosen sets.
The contracted graphs only lose the chosen vertices, so a core
separator remains one after the chosen vertices are removed from it.
Each core $(i,j)$ is directly anchored by exactly one permutation
type, namely $X_{i-j}$, so its left-anchor count and its left
deficit both decrease by $h$, and activity is preserved.
If $\min_i a_i\le\min_s p_s$, taking $h=\min_i a_i$ gives
a recoverable target. Thus we may take $h=\min_s p_s$.
\end{proof}

\begin{lemma}\label{r3.3}\label{r3.4}
A balanced bipartite graph on $m+m$ vertices, missing at most
$e$ edges at each vertex, has a Hamilton path between any
prescribed opposite endpoints if $m-1\ge4e$ and $m-e\ge2$.

Let a monochromatic cycle $C_0$ of half-length $h_0$ contain
an edge of every selected, pairwise disjoint core $D$. If it
uses $u_D,v_D$ core vertices and the missing-degree bound is
$e_D$, put $T_D=\min\{|A_D|-u_D,|B_D|-v_D\}$. For integers
$0\le t_D\le T_D$ with every positive $t_D\ge4e_D$, the cycle
extends to half-length $h_0+\sum_Dt_D$. Complete cores therefore
give every integer half-length up to $h_0+\sum_DT_D$.
\end{lemma}
\begin{proof}
Order the left vertices $a_1=x,\ldots,a_m$, with
$a_m\in N(y)\setminus\{x\}$. Match consecutive pairs
$(a_i,a_{i+1})$ to right vertices other than $y$, requiring
adjacency to both entries. Both sides of this auxiliary graph
have minimum degree at least $m-1-2e\ge(m-1)/2$: a missing
left vertex affects at most two consecutive pairs. Hall's
theorem gives the path $a_1b_1\cdots a_{m-1}b_{m-1}a_my$.
For each positive increment, replace one edge in $D$ by this
path on its endpoints and $t_D$ unused vertices of each shore.
The replacements are disjoint; zero increments leave edges
unchanged. The complete-core capacities are integer intervals.
\end{proof}

For each colour, contract its three large cores, retaining corner
edges and exceptional incidences. An exceptional vertex adjacent
to a core is a direct \emph{anchor}. Write $L_K,R_K$ for the
left and right anchor sets of $K$. A \emph{core separator} is
a set of exceptional vertices whose deletion separates the cores.

\begin{lemma}\label{r3.5}
A target exists if a contracted colour graph has a cycle through
two cores, or if a family of vertex-disjoint exceptional paths
with endpoints anchored at one core $A_i\times B_j$ contains
at least $a_i$ left and $b_j$ right vertices. Anchor singletons
are allowed. Hence no monochromatic cycle in a counterexample
visits two cores.
\end{lemma}
\begin{proof}
Lift core visits to paths with distinct entries. For a one-core
path family, join the paths cyclically inside that core. Its
shores have size at least $n-r>2r$, so entries are available
and the resulting half-length satisfies $h_0\le2r<n$.
Two visited cores each use at most two vertices per shore, giving
$T_1+T_2\ge2(n-r-2)\ge n$. For one core with exceptional
counts $\ell,s$, the remaining capacity is
\[
 T=n-h_0+\min\{\ell-a_i,s-b_j\}\ge n-h_0.
\]
Apply \cref{r3.4}. Finally, contracting a simple cycle through
several cores gives an even multigraph with degree-two noncore
vertices. Delete the core loops and the exceptional segments
returning to the same core. The remaining segments form a closed
walk through at least two distinct cores, which contains a cycle
through distinct cores, allowing parallel paths between two cores.
Its exceptional paths have disjoint interiors.
\end{proof}

\begin{lemma}\label{core:rules}\label{r3.19}
Let $T$ be an independent set of three core vertices in a graph
with no cycle through two cores.
\begin{enumerate}
\item The core--anchor incidence graph is a forest. If $d(u)$
is an anchor's number of core neighbours, then
$\sum_u(d(u)-1)_+\le2$.
\item For distinct $D,Q\in T$ and disjoint anchor sets
$P\subseteq N(Q)\setminus T$, $Y\subseteq N(D)\setminus T$,
the edges between $P$ and $Y$ form a possibly empty star.
\item In a normalised Latin counterexample, mark $h$ cores
with $b(K)=1$, $R_K=\varnothing$ and
$a(K)\le\max\{2,|L_K|\}$. There is a core separator with
\[
 |Z_c|\le2,\qquad |Z_c\cap R_*|\le(2-h)_+.
\]
Thus all colours need at most $6-\sum_c\min\{2,h_c\}$ right
memberships.
\end{enumerate}
\end{lemma}
\begin{proof}
Each block contains at most one core: two vertices of a
nontrivial block lie on a common cycle, and no bridge joins
two cores. The incidence graph is therefore a forest, whose
edge count, after discarding isolated anchors, gives (i).
For (ii), two independent cross-anchor edges give the cycle
$Q-x_1-y_1-D-y_2-x_2-Q$.

A cycle at a marked core gives an anchored path with at least
two left vertices and one right vertex. Adding unused left
anchors gives absorption by \cref{r3.5}. Every edge at that
core is therefore a bridge, and all its neighbours are left.
Root each core-containing block--cutpoint tree at an unmarked
core when possible, representing noncut cores by leaves.
For each nonroot core, choose the first exceptional cutvertex
towards the root. These choices separate the cores; a marked
core chooses a direct left neighbour. If $k$ trees contain
cores and $k_0$ contain an unmarked core, then
\[
 |Z_c|\le3-k\le2,\qquad
 |Z_c\cap R_*|\le3-h-k_0\le(2-h)_+,
\]
since $k_0\ge\ind{h<3}$. Repeated choices only improve the bounds.
\end{proof}

\begin{lemma}\label{r3.6}
If three cores in distinct rows have left-anchor sets covering
$L_*$, some core has fewer right anchors than its right deficit.
If only $D$ does, at least $a_D+2$ left anchors are exclusive
to $D$. For a nonempty same-colour family $S$ of such cores,
\[
 |L_D|\ge a_D+2,\qquad
 \left|\bigcup_{D\in S}L_D\right|
 \ge1+\sum_{D\in S}(a_D+1).
\]
Every column is a covering triple, as is any colour used by
every left vertex. The assertions transpose.
\end{lemma}
\begin{proof}
A right-sufficient core has at most $a_D-1$ left anchors by
\cref{r3.5}. Since the deficits sum to $r_L$, not all three
cores can be right-sufficient. If only $D$ is deficient,
subtract the other two anchor sets to obtain its exclusive
anchors. By \cref{core:rules}(i),
$\sum_{D\in S}|L_D|\le|S|+|\bigcup_{D\in S}L_D|-1$,
which gives the displayed bounds.
\end{proof}

\begin{lemma}\label{r3.11}
Let $P$ be a left permutation type with principal core $Q$ in
colour $c$. Suppose two connections on distinct exceptional
vertices form a tree on the three colour-$c$ cores. A connection
is an ordinary vertex's direct pair, or a right permutation
vertex directly anchored at $D\ne Q$ with two $P$ neighbours.
Then $N_c(y)\cap P=\varnothing$ if $y$ is not a connection
vertex but is directly anchored at another core, or if $y$ is
a repeated connection vertex whose pair omits $Q$. No other
vertex repeats $c$.
\end{lemma}
\begin{proof}
Tree paths lift using distinct $P$ vertices. A path ending at
$Q$ uses at most one indirect connection and can avoid a
specified $x\in P$. An edge $xy$ then closes $Q-x-y-D$ with
this path to give a two-core cycle. For a connection vertex
omitting $Q$, use the other connection. Any further repeated
vertex similarly closes its direct pair with the tree.
Repeated connection vertices are not in $P$, so no exceptional
vertex is reused.
\end{proof}

\begin{lemma}\label{r3.12}
For a core $D$, let $L_D,R_D$ be its anchor sets and put
$d_D=|R_D|$. Every exceptional left vertex $x$ satisfies
\[
 |N_c(x)\cap R_D|\ge2,\quad a_D\le|L_D\cup\{x\}|
 \quad\Longrightarrow\quad d_D<b_D.
\]
If $P$ anchors another core and $Y\subseteq R_D$, two nonempty
$P$-neighbourhoods in $Y$ must be the same singleton.
If $a_D=1,b_D\le2$, at most one vertex of $Y$ meets $P$,
and every exceptional left vertex has at most one neighbour in
$Y$. All assertions transpose.
\end{lemma}
\begin{proof}
For the displayed implication, take two neighbours $y,z\in R_D$ and use
the path $y-x-z$ with unused anchor singletons. If $d_D\ge b_D$,
this gives absorption. The singleton assertion is
\cref{core:rules}(ii). In the unit case, two neighbours in $Y$
would give $d_D\ge2\ge b_D$, contradicting the displayed implication.
The star rule then also excludes two vertices of $Y$ meeting $P$.
\end{proof}

A core is \emph{active} if its direct left-anchor count meets
its left deficit. A permutation type is \emph{fully active}
if all three principal cores are active. For an active core
put $k=b-d$, where $b$ is its right deficit and $d$ its direct
right-anchor count. By \cref{r3.5}, if $k=1$, one further
right vertex with two left-anchor neighbours gives absorption.
If $k=2$, two further vertices with two and three such neighbours
suffice. Starting with anchor singletons, each insertion joins
two paths; three choices ensure distinct paths at the second
insertion.

\begin{lemma}\label{r3.23}
Suppose $P$ is fully active and every $b_j\le2$, and put
$p=|P|$. Let $\beta_c,d_c$ be its principal colour-$c$ core's
right deficit and right-anchor count, and set
$k_c=\beta_c-d_c\in\{1,2\}$. There are integers
$1\le t_c\le k_c$ and $Y\subseteq R_*$ avoiding all principal
cores such that
\[
 |Y|\ge T:=\sum_ct_c\ge p,\qquad |N_c(y)\cap P|\le t_c.
\]
In particular, $p\le r_R-\sum_cd_c$; if $r_R\le4$, then $p\le3$.
\end{lemma}
\begin{proof}
For $k_c=1$, insertion bounds nonanchor degrees by one.
For $k_c=2$, either all degrees are at most two, or at most
one exceeds one. Take $t_c=2$ in the first case and $t_c=1$
in the second, deleting the latter vertex. Delete the direct
anchors as well. At most $\beta_c-t_c$ vertices are deleted
for each colour, so $|Y|\ge r_R-\sum_c(\beta_c-t_c)=T$.
Summing the colour degrees of any $y\in Y$ gives
$p\le T\le\sum_ck_c$.

If $r_R\le4$ and $p\ge4$, equality forces
$p=|Y|=T=r_R=4$, with thresholds $(2,1,1)$ attained at every
vertex. In the threshold-two colour $c_0$, \cref{r3.12}
makes permutation types distinct and forbids a repeated-$c_0$
vertex sharing a nonprincipal core with a permutation.
There are no constants and at most one ordinary vertex per
repeated colour. Hence exactly two incompatible permutations
and two ordinary vertices occur. In every colour, the two
permutations' direct incidences, together with one $P$ edge
from each, form a tree. Either ordinary repeated pair closes
it to a two-core cycle.
\end{proof}

\begin{lemma}\label{r3.26}
If $r_R\le4$, a left deficit is one, and a repeated left
vertex exists, every fully active left permutation type has
size at most two.
\end{lemma}
\begin{proof}
Otherwise \cref{r3.23} gives $p=3$. If $T=3$, every
$y\in Y$ uses all three colours towards $P$. In a row of
deficit one, its coordinate belongs to one of the two
nonprincipal cores, since $Y$ avoids the principal cores.
The unit-row case of \cref{r3.12} permits at most one such
vertex meeting $P$ at each of these cores, giving $|Y|\le2$.
Thus $T=|Y|=r_R=4$, with thresholds $(2,1,1)$, say doubled
in $c_0$, and $Y=R_*$. No constant vector avoids all principal
cores. Every right vertex has positive $c_0$ degree into $P$,
since its other two degrees are at most one. If two right
permutations had the same type, they would anchor the same
nonprincipal $c_0$ core; the star rule in \cref{r3.12} would
make both $c_0$ neighbourhoods the same singleton. Their total
degree is three, so both would use all three colours towards
$P$, contrary to the unit-row case of \cref{r3.12} at their
common coordinate in that row. Hence their types are distinct.

Let $v$ count these permutations, and let $\delta$ indicate
the presence of an ordinary vertex repeating $c_0$.
Ordinary repeated colours are distinct, so $v\ge2-\delta$.
If $\delta=1$, its shared nonprincipal core bounds each
permutation's $c_0$ degree by one. Each permutation therefore
uses at least $2+\delta$ colours towards $P$. Add the
corresponding $P$ edges to the direct forests. Each colour
allows at most two connections: as the permutation types are
distinct, three connections on its three cores would close a
two-core cycle, shortcutting through a common vertex of $P$ if
necessary. The $4-v$ ordinary right vertices and a repeated left
vertex also contribute connections. Hence
\[
 6\ge1+(4-v)+(2+\delta)v
 \ge5+(1+\delta)(2-\delta)=7.
\]
\end{proof}

\begin{theorem}\label{r3.29}
If $r_R=3$ and every $a_i\le2$, a target exists.
\end{theorem}
\begin{proof}
Suppose otherwise. Positivity gives $b_0=b_1=b_2=1$.
Let $M_c$ count exceptional vertices missing
$c$, and let $f_c$ be the repeated weight in colour $c$: one
per ordinary vertex repeating $c$, two per constant-$c$ vertex.
Put $t_c=M_c+1-f_c$. The direct forests give $f_c\le2$.
Every vector has equal missing and repeated weight, so
$\sum_ct_c=3$.

Suppose first that $t_c<0$. Then $f_c=2,M_c=0$, and the direct
colour-$c$ graph $S$ is a spanning three-core tree. By
\cref{r3.6}, one or two cores have right anchors. If only one does,
the two other ordinary repeated colours, each occurring at most
once, cannot both occur: they would make a whole column
right-sufficient. Thus at least two right vertices are
permutations, with value $c$ in the row of that core. Of the two
permutations with this value there, one is a class vector,
excluded by normalisation, so these right vertices are equal.
If all three formed one permutation type $M$, it would be
fully active after transposition. Transposed \cref{r3.23}
gives a left set $Y$, avoiding all its principal cores, and
thresholds $1\le s_d\le2$ with
$|Y|=m\ge T=\sum_ds_d$ and $|N_d(y)\cap M|\le s_d$.
As $M_c=0$, every $y\in Y$ has a colour-$c$ coordinate;
its colour-$c$ core is not the principal core of $M$.
A colour-$c$ edge from $y$ to $M$ would close the path in
$S$ to a two-core cycle. Thus $|N_c(y)\cap M|=0$.
The other two degrees sum to three and are each at most two,
so both are positive and their thresholds sum to at least
three. Hence $T\ge4$ and $m\ge4$.
For either other colour $d$, every $y\in Y$ using a
nonprincipal coordinate of colour $d$ also has an actual
colour-$d$ edge to $M$. The transposed unit-row case of
\cref{r3.12}, using $b_j=1$ and $a_i\le2$, permits at most
one such $y$ at each coordinate. There are only four such
coordinates. Every $y$ uses at least two, except the ordinary
vertices repeating $c$, which use one; there are at most
$f_c=2$ of these. Thus $4\ge2m-2\ge6$, a contradiction.
Using the symmetries $(i,j,c)\mapsto(i+u,j+v,c+u+v)$ and
$(i,j,c)\mapsto(-i,-j,-c)$, we may therefore take $c=0$ and right
vectors $Y_0,Y_0,(0,1,1)$. Applying \cref{r3.6} to column zero gives
$a_2+2\ge3$ left anchors $G$ of $(2,0)$. Their zero edges to $R_*$ are forbidden
by $S$, and each right vertex has at most one colour-two
neighbour in $G$. On three such left vertices, colour one has
a path through all three right vertices, since each right
vertex has at least two colour-one neighbours in that triple.
These anchor $(2,2)$,
and the two internal left vertices meet $a_2\le2$, giving absorption.

If two cores have right anchors, the third $D$ has $a_D+2$
exclusive left anchors $Z$. Its ordinary members repeat each
other colour at most once. If $a_D=1$, both cannot occur,
since they would activate all three unit-row cores, contrary to
transposed \cref{r3.6}. Thus $Z$
contains two equal permutation vertices of a fully active type
$P$. The other colour degrees into $P$ are at most one and
$S$ forbids $c$, so $|P|=2$, with both other colours used at
every right vertex. Take $P=X_0,c=0,D=(0,0)$. The ordinary
vertices are $(0,d,d)$ for $d\in I\subseteq\{1,2\}$, with
$a_0\le|I|$. Each joins $P_d$ to $(0,d)$, so every right
zeroth coordinate avoids $I\cup\{0\}$. Hence $|I|=a_0=1$,
and all these coordinates equal the remaining colour $e$.
Two right vertices share an $e$ neighbour in $P$, absorbing
$(0,e)$. We have proved $0\le t_c\le3$.

Let $\nu_c,\kappa_c$ be the corner cycle rank and component
count, including isolated vertices, and let $I_c$ count direct
incidences. Put $N=r_L+3$. Since the corner graphs have
$3r_L$ edges in total, $\sum_c\nu_c=\sum_c\kappa_c-9$.
Contract corner components and merge parallel incidences,
with excess $\Delta_c$. The quotient forest gives
$I_c-\Delta_c\le\kappa_c+2$. As $\sum_cI_c=3N$,
$\sum_c\nu_c\ge3r_L-6-\sum_c\Delta_c$.
Unit cores have no excess: two incidences from one corner
component give a path meeting both unit deficits. At a
deficit-two core, two left anchors exclude right anchors by
singleton absorption. They also lie in distinct corner
components, since a path between them already has two left
vertices and one right vertex. Otherwise the excess is at most
$h-\ind{\lambda=0,h>0}$ for $\lambda$ left and $h$ right
anchors. Let $z_c$ count such cores with right but no left
anchors. Each of the $r_L-3$ deficit-two rows contributes
three right incidences. Thus
$\sum_c\Delta_c\le3(r_L-3)-\sum_cz_c$, giving
\[
 \sum_c\nu_c\ge3+\sum_cz_c.
\]

If $\nu_c>0$, there is a unique cyclic corner component $J$,
since only three right vertices exist. At each deficit-two core
without left anchors, retain one right incidence per corner
component. A repeated component contains at least two right
vertices, as does $J$, so it must be $J$. If all three right
anchors lay in $J$, a cycle in $J$ would give a path through
all three: this is immediate for a six-cycle, while for a
four-cycle one extends a shortest path from the third right
vertex to that cycle. This path would absorb the core.
Thus at most one incidence is removed per core.
After $d\le z_c$ deletions, the direct
forest $F$ has $t_c+d+2$ components. A corner edge whose
endpoints share a core neighbour in $F$ is a bridge, since an
alternative path gives absorption. Delete these bridges.

Every remaining corner component meets each $F$ component
at most once. Indeed, take a shortest return path with interior
outside that $F$ component. Its union with the $F$ path gives
a two-core cycle or an absorbing one-core path. A single edge
was already deleted. The only short exception is a
right--left--right path at a deficit-two core: add an unused
left anchor; if its sole left anchor is the middle vertex,
both edges were deleted; with no left anchor, the retained
incidence rule excludes the path.

Thus the remaining cyclic component has $l+s\le t_c+d+2$
vertices. If $d=0$, this is at most five, giving
$\nu_c\le l+s-3$. If $d>0$, two left vertices cannot both
see all three right vertices: at a core with a deleted
incidence, $y_1-x_1-y_3-x_2-y_2$ would absorb. Hence
$e(J)\le2l+s-2$: for $s=3$, at most one left vertex has
degree three, and for $s=2$ all degrees are at most two.
As $J$ is cyclic, $s\ge2$, so
$\nu_c=e(J)-l-s+1\le l-1\le l+s-3$, giving the same bound.
Therefore
$\nu_c\le(t_c-1)_++z_c$, also when $\nu_c=0$.
Since $\sum_ct_c=3$, summing gives at most
$2+\sum_cz_c$, contrary to the lower bound $\sum_c\nu_c\ge3+\sum_cz_c$
obtained above.
\end{proof}

\begin{lemma}\label{r3.30}\label{r3.32}
A counterexample has $r\ge8$. If $r_L=3$ and two left
permutation types occur, then $r_R\ge6$.
\end{lemma}
\begin{proof}
If $r\le7$, transpose to get $r_R=3,r_L\le4$, contrary to
\cref{r3.29}. For the second assertion, suppose $r_R\le5$.
Transposed \cref{r3.29} forces $\boldsymbol b=(3,1,1)$.
Take $P=X_1,Q=X_2$. Their unit-column cores are active, so
the right vectors are $(0,2,1),(0,1,1),(0,2,2)$; the remaining
candidate $(0,1,2)$ is the class vector of $B_0$.
The last two occur at most once, leaving at least three
permutations $M=(0,2,1)$. Their right-anchor counts meet
all deficits, so \cref{r3.12} bounds each exceptional left
colour degree into $M$ by one. Thus $|M|=3$, and every left
vertex uses every colour. Both ordinary vectors occur;
$(0,1,1)$ joins $P_1$ to the colour-one core of $M$, and a
$P$--$M$ edge closes the path.
\end{proof}

In a fixed colour, a right vertex has an \emph{effective
incidence} with a core if it is directly anchored there or has
two neighbours in one left permutation type anchored there.
Given core separators $Z_c$, a vertex is \emph{interior} if it
lies in $Z=\bigcup_cZ_c$ and \emph{exterior} otherwise.

\begin{lemma}\label{r3.36}
If no right vertex belongs to two core separators $Z_c$, a
fully active permutation type has at most two vertices outside
$Z=\bigcup_cZ_c$.
\end{lemma}
\begin{proof}
Suppose $P_0=P\setminus Z$ has size at least three, and put
$k_c=\beta_c-d_c\ge1$. The right set $U$ avoiding the
principal cores has size at least
$r_R-\sum_cd_c=\sum_ck_c$. No constant vertex lies in $U$.
Let $r_c\le1$ count its ordinary
vertices repeating $c$. If a permutation $y\in U$ has a
colour-$c$ neighbour $x\in P_0$, its direct colour-$c$ core
differs from the principal core of $P$. The path through
$x,y$ between these cores would survive deletion of $Z_c$
unless $y\in Z_c$. Since no right vertex lies in two
separators, all edges from $y$ to $P_0$ have this one colour
$c$. In particular, $y\in Z_c\cap U$ has at least three
colour-$c$ neighbours in $P_0$ and an effective incidence
with the principal core of $P$. Let $m_c$ count these permutations.
Their effective connections and repeated pairs
form a three-core forest, so $r_c+m_c\le2$. For $k_c\le2$,
insertion gives $m_c\le k_c-1$; for $k_c\ge3$, the forest
bound is strict. Thus
\[
 \sum_ck_c\le|U|=\sum_c(r_c+m_c)\le\sum_ck_c,
\]
forcing $k_c\le2,r_c=1,m_c=k_c-1$. Choose an ordinary
$y\in U$. If it has two neighbours of one colour in $P_0$,
insert it first, then the completely joined vertices, giving
absorption. Each colour degree is at most one, and its single
colour is forbidden by separation. It cannot see three vertices
of $P_0$.
\end{proof}

\begin{lemma}\label{r3.35}
The union of any core separators with $|Z_c|\le2$ contains
all repeated vertices and covers one exceptional shore. After
transposition,
\[
 R_*\subseteq Z:=\bigcup_cZ_c,\qquad
 3\le r_R\le6,\qquad w\le6-r_R,
\]
where $w$ counts repeated left vertices.
\end{lemma}
\begin{proof}
Separators exist by \cref{core:rules}(iii), and repeated
pairs force their vertices into the corresponding separator.
If both shores have exterior vertices, their types are unique
and compatible. Let their sizes be $p,q$, and let $z_L,z_R$
be the interior counts. Put $x=|\{i:a_i>p\}|$ and
$y=|\{j:b_j>q\}|$. Avoiding absorption and counting positive
deficits give
\[
 x+y\ge3,\quad z_L\ge3+(x-1)p,\quad z_R\ge3+(y-1)q,
 \quad (x-1)p+(y-1)q\le0.
\]
After transposition, $x=0,y=3,q=1,z_R\ge5,z_L\le1$.
Also $\max a_i\le p+z_L-2<p$. If $p=2$, then
$r_L=3,r_R=6,\boldsymbol b=(2,2,2)$, excluded by transposed
\cref{r3.29}. For $p\ge3$, $|Z|=6$ contradicts
\cref{r3.36}. Otherwise $z_L=0,z_R=5,r_R=6,b_j=2$.
The exterior right vertex anchors each principal core once;
any further direct anchor or two-neighbour right vertex gives
absorption. Thus $p=3$, and every interior right vertex uses
all colours towards $P$. Its vector uses at least two colours,
since a constant would give a forbidden principal anchor.
These require two separator memberships each, so five vertices
need ten positions, exceeding six.
\end{proof}

\emph{Full separation} means six distinct separator positions
consisting exactly of $R_*$ and the repeated left vertices.
Left permutation vertices are then exterior, and each right
vertex has one assigned colour.

\subsection{The six-position budget}\label{sec:4}
We say that \emph{one-sided coverage} holds if there are core
separators $Z_c$ with $|Z_c|\le2$ whose union contains $R_*$.

\begin{theorem}\label{r4.2}
One-sided coverage and $r_R\ge5$ force a target.
\end{theorem}
\begin{proof}
Apply \cref{r3.2}, leaving types $P,Q,U=X_0$ of sizes
$p\ge q\ge0$, with $U=\varnothing$. Minimise the right
membership count $\sigma$ among separators covering $R_*$.
Let $\ell$ count interior left permutation vertices, and let
$E=\sum_v(m(v)-1)_+$, where $m(v)$ is the number of separators
containing $v$. Then
\[
 r_R+w+\ell+E\le6,\qquad w+\ell+E\le1.
\]
\cref{r3.32} gives $r\ge8,p\ge2$; otherwise
$r_L=3,p=q=w=1,r_R\le5$.

A \emph{three-arm star} is a right vertex of some $Z_c$ joined
in colour $c$ to all three colour-$c$ cores, each directly or
through one left vertex; these three paths are its arms.
A three-arm star excludes another right vertex with two effective
incidences, or another vertex repeating its colour: the
corresponding tree paths close to a two-core cycle. If its $P$ arm
has two choices and its $Q$ arm ends at $v$, deleting the centre
and $v$ separates the cores: a shortest remaining core path
contains at most one vertex of $P$, as every vertex of $P$ is
adjacent to $P_c$, so it closes with an unused $P$ arm.

Consider an exchange replacing one separator by another of
size at most two, introducing only left members and removing
exactly $k\ge1$ right memberships. It decreases $\sigma$ by $k$
and retains at least $r_R-k$ distinct right vertices in the
union. Minimality excludes continued right coverage, while
\cref{r3.35} forces left coverage. Hence $r_L+r_R-k\le6$,
so $r\le6+k$; in particular, $k=1$ is impossible.

Suppose $y\in Z_c$ has at most one effective incidence.
If $w=0$, replace $y$ by one singleton neighbour unless its
direct core and neighbouring types cover all three cores.
In the latter case it directly meets $U_c$ and has unique
neighbours $u\in P,v\in Q$. In the contracted colour-$c$
graph, $y$ has precisely the neighbours $U_c,u,v$.
The set $\{u,v\}$ itself separates the cores: after its deletion
$y$ is a leaf, so any remaining path between cores avoids $y$.
The subpath between two consecutive cores, together with the
corresponding arms of the original star, would then form a
two-core cycle, contrary to \cref{r3.5}.
Replace the whole separator $Z_c$ by $\{u,v\}$, and put
$k=|Z_c\cap R_*|\in\{1,2\}$. The case $k=1$ is already
excluded. For $k=2$, the exchange bound gives $r\le8$;
as $r\ge8$, $r_L\ge3$ and $r_R\ge5$, this forces
$r=8,r_L=3,r_R=5$. Both $P,Q$ are nonempty, contradicting
\cref{r3.32}.
If $w=1$, then $r_R=5,\ell=E=0$. Such a $y$ must be a
permutation compatible with $P$, since incompatible $P$ edges
have its assigned colour and give a second incidence.
Its $Q$ edges give $q\le1$, so
$\max a_i\le p+q-1\le p$; compatibility forces all
$b_j\ge2$, a contradiction. Every right membership is
therefore at least doubly effective, and two in one colour
cannot share two incidences.

Let $e_c$ count right members of $Z_c$, let $\mathcal D$ be the
set of cores with no effective incidence from a right member of
the separator of their colour, and let $J$ be their columns.
Since $\sigma\ge5$ and $e_c\le2$, every $e_c\ge1$.
These members cover at least $e_c+1$ cores, so
$|\mathcal D|\le6-\sigma$. Every column has an active core;
outside $J$, effective coverage forces its deficit to be at
least two. Hence $r_R\ge2(3-|J|)+|J|$, and
\[
 6-r_R\le|J|\le|\mathcal D|\le6-\sigma\le6-r_R.
\]
All are equal, yielding
\[
 b_j=2-\ind{j\in J},\quad E=0,\quad w+\ell\le1.
\]
Every right vertex has exactly two effective incidences.
An active core has no direct right anchor outside its
separator, since this fills its unit deficit or supplements
an effective incidence. Transposed \cref{r3.29} gives
$r_L\ge4$. Both $P,Q$ have exterior vertices; otherwise a
sole exterior type has at least $r_L-1\ge3$ vertices and
exceeds every $a_i$, contrary to \cref{r3.36}.

A right permutation compatible with a type of size $h$
anchors it outside its separator in two rows, both inactive.
For the other size $v$, this gives
\[
 v\ge h+3-w,\qquad \max a_i\le v+w-2<v,
 \qquad v-\ell\ge3,
\]
contrary to \cref{r3.36}. A remaining right permutation is
compatible with $U$. Its three-arm star leaves one right
member in $Z_c$ and forbids repeated-$c$ vertices. Thus
$r_R=5,w=0$, the other separators are right pairs, and all
left memberships lie in $Z_c$. Hence an edge from $y$ to an
interior vertex $x_0$ of a type $V$, of colour $c'\ne c$, would
give the path $V_{c'}-x_0-y-U_{c'}$ surviving the deletion of
$Z_{c'}$, and the same holds for exterior vertices. Thus $y$ is
joined to all of $P\cup Q$ in colour $c$. As $p\ge2$, this gives
an effective incidence with $P_c$, so exactly two effective
incidences force $q\le1$, and $q\ge1$ since $Q$ has an exterior
vertex. Write $Q=\{v\}$ and exchange $Z_c$ for
$\{y,v\}$. The active type $P$ has $p=r_L-1\ge3$ exterior
vertices, again a contradiction.

All right vertices are therefore ordinary, assigned their
repeated colours. Two colours form three-core trees; at most
one core lacks a repeating right anchor, and every core has
at least $b_j-1$ such anchors. By \cref{r3.11}, a tree
vertex $z$ omitting type $V=X_v$ has no repeated-colour edge
to $V$. Write $c,d$ for its repeated and single colours.
Its single-colour core lies in the omitted row $s=-c-v$,
whereas $V_d$ lies in row $-d-v\ne s$. Since $z$ belongs
only to its repeated-colour separator $Z_c$, a colour-$d$
edge from $z$ to an exterior vertex of $V$ would connect these
two distinct cores after deleting $Z_d$.
Its single-colour neighbours therefore lie among the $\ell$
interior vertices, leaving at least $|V|-\ell$ neighbours in
its missing colour. If $p-\ell\le1$, then
$p=q=2,\ell=1,w=0,r_R=5$; two active unit-column cores
lack right anchors, a contradiction. Hence $p-\ell\ge2$.

If $P$ is fully active, some tree vertex omits it; otherwise
two direct anchors give absorption. Its missing-colour path
and the $b_j-1$ direct anchors also give absorption. Thus
$a_{i_0}>p$ for some row, and
\[
 q\ge3-w,\quad a_i\le q+w-2\le q-1\ (i\ne i_0),
 \quad q-\ell\ge2.
\]
Both types are active outside $i_0$. In each tree colour $c$,
choose a vertex omitting $V=X_v\in\{P,Q\}$, with single
colour $d$, missing colour $e$, and omitted row $s$.
Its missing-colour path forces $V_e$ into row $i_0$, so
$v=-e-i_0$ and $s=-c-v\ne i_0$. Its single-colour core
has its own anchor and $b_j-1$ others, so is inactive and
belongs to $U=X_0$, giving $s=-d$. Since $c+d+e=0$,
we obtain $c=-i_0$ for two distinct colours, a contradiction.
\end{proof}

\subsection{The Latin-core theorem}\label{sec:5}
\begin{theorem}\label{r5.7}
Every homogeneous Latin model with positive deficits and
$n\ge3r+6$ contains a target.
\end{theorem}
\begin{proof}
Normalise a counterexample. By \cref{r3.35,r4.2}, choose
separators covering the right side with $3\le R=r_R\le4$,
$R+w\le6$ and $r\ge8$. We first establish three local
facts, also used after transposition.

\emph{(A) Under full separation, a fully active type has size
at most one.} Otherwise \cref{r3.36} gives
$P=\{x_1,x_2\}$. Since the positive right deficits sum to
$R\le4$, at least two are one. A right permutation compatible
with $P$ directly anchors every principal core of $P$, including
an active core with right deficit one; its singleton, together
with enough left-anchor singletons, gives absorption by
\cref{r3.5}. Thus every right permutation is incompatible
with $P$. If such a vertex $y$ has a colour-$c$ neighbour
$x\in P$, its direct colour-$c$ core differs from $P_c$.
As $x$ lies outside every separator, the path from $P_c$
through $x,y$ to that core forces $y\in Z_c$. Hence every
right permutation uses its unique assigned colour towards $P$.

Each $Z_c$ has two members. A repeated member supplies its
direct pair, while an assigned right permutation supplies a
connection from its direct core to $P_c$ through either $x_i$.
Coincident connections give a two-core cycle, using distinct
$x_i$ when needed. A constant member's three-arm star, together
with the other connection, also gives such a cycle. Thus all
repeated members are ordinary and the two connections form a
tree. There are $u=6-R=w$ unit principal deficits. For each
such colour $c$, a direct right anchor at $P_c$, or an assigned
right permutation with its two $P$ neighbours, gives absorption.
Consequently any right member of $Z_c$ is ordinary and its
repeated pair omits $P_c$. Two such members would have the same
pair. Each unit-colour separator therefore needs a left member;
as there are exactly $u$ left members altogether, it has exactly
one left member $L_c$ and one right member $R_c$. The pair of
$L_c$ meets $P_c$, since otherwise it coincides with that of
$R_c$.

Write $P=X_v$, and let $d,e$ be the single and missing colours
of $R_c$. Its repeated pair omits $P_c$, whose row is
$s=-c-v$, so its single colour occurs in row $s$. The row of
$P_d$ is $-d-v\ne s$. Thus separation forbids colour $d$
between $R_c$ and $P$, while \cref{r3.11} forbids colour $c$.
Both edges must have colour $e$. If the principal right deficit
of colour $e$ were one, the path $x_1-R_c-x_2$, with unused
left-anchor singletons, would give absorption. Hence $e$ is
nonunit, excluding $R=3$.

For $R=4$, take $P=X_0$ and nonunit colour zero. Then
$R_1=(1,1,2)$ and $R_2=(2,1,2)$. Both members of $Z_0$
are right vertices. An ordinary member whose repeated pair
omits $P_0$ has only its missing colour towards $P$, by the
same coordinate argument; that colour is unit, giving
absorption. An ordinary member whose pair meets $P_0$ gives
one direct right anchor there. The other member either gives
a second direct anchor, or is a permutation whose two $P$
neighbours fill the remaining right deficit. Thus $Z_0$
consists of two distinct permutations
$Y_1=(1,0,2)$ and $Y_2=(2,1,0)$.

If another left type $Q$ occurs, take $v\in Q$. Exactly one
of $Y_1,Y_2$, say $y$, is compatible with $Q$; call the other
$z$. Separation forces $zv$ to have colour zero, so
$P_0-x_1-z-v-Q_0-y-x_2-P_0$ is a two-core cycle.
Hence $r_L=4$. In colour one, the two nonprincipal cores
are in rows zero and one and have respectively two and three
right anchors, enough to meet their right deficits. Their
left-anchor counts sum to at least one, supplied by $L_1$,
so \cref{r3.5} gives $a_0+a_1\ge3$. Similarly,
$a_0+a_2\ge3$. Positivity and $\sum_i a_i=4$ now give
$\boldsymbol a=\boldsymbol b=(2,1,1)$.

In column zero, the colour-one and colour-two cores each have
three right anchors, whereas the colour-zero core has none.
Applying \cref{r3.6} gives $a_0+2=4$ exclusive left anchors
at the latter core. Thus $L_c=(0,c,c)$, and all right vertices
are joined to $P$ in colour zero. No $L_cy$ has colour zero:
for $y'\ne y$, the path $L_c-y-x_1-y'-x_2$ is anchored at
$(0,0)$ and meets its deficits $(2,2)$. Put
$D_c=(0,c)$ and $E_c=(c,0)$. Neither $L_1R_2$ nor $L_1Y_2$
has colour one, since either closes the two-core cycle
$L_1-D_1-R_1-E_1-z-L_1$, with $z=R_2$ or $Y_2$.
Both edges therefore have colour two. The $D_2$-anchored path
$R_2-L_1-Y_2$, together with singleton $L_2$, has two left
and two right vertices and absorbs the deficits $(2,1)$.

\emph{(B) Two types have at most $u=6-R$ active unit-column
cores.} Suppose that there are $m>u$ such cores. Each has no
direct right anchor, since it is active and has right deficit
one, and its left-anchor count meets its left deficit. They
may all be marked in \cref{core:rules}(iii). Each colour has
at most two marked cores, so the resulting separators have at
most $6-m<R$ right memberships and omit a right vertex.
Their union covers the left shore by \cref{r3.35}.
Transposed \cref{r4.2} gives $r_L\le4$; together with
$r_L+R\ge8$ and $R\le4$, this forces $r_L=R=4$.

At most one right permutation type occurs. Indeed, two types
directly right-anchor two distinct cores in every unit column,
leaving at most one active core per such column, contrary to
$m>u$. The separator union contains all four left vertices
and every repeated right vertex, so at least two right vertices
are permutations of the same type $M$. After transposition,
all left deficits are at most two, so $M$ is fully active.
Transposed \cref{r3.23} gives $|M|\le3$. If equality holds,
the fourth right vertex is repeated, and the unit left deficit
allows transposed \cref{r3.26} to give $|M|\le2$ instead.
Thus $|M|=2$, and the four left vertices together with the two
repeated right vertices occupy six distinct separator positions.
This is transposed full separation, contradicting (A).

A fully active type already contributes $u$ active unit-column
cores, so another type has none. If two types are active in two
common rows, their four active cores exclude $R=3$, when all
columns are unit. For $R=4$, at most two of these four cores
are unit; hence the nonunit column contains one of the specified
cores of each type.

\emph{(C) For $\boldsymbol b=(2,1,1)$, if $X_1,X_2$ are
active in rows $1,2$, then $a_1+a_2>p$, where $p$ is their
larger size.} Suppose not, and take $P=X_1$, interchanging
indices one and two if necessary. The active colour-zero cores
in rows two and one have right deficit one, so neither row
coordinate of a right vector can be zero. The active cores
$(X_1)_1$ and $(X_2)_2$ have right deficit two, so each has
at most one direct right anchor. Hence at most one right vector
takes value one in row one, and at most one takes value two in
row two. Let $M$ consist of all right vertices taking values
$(2,1)$ in these two rows. Then $|M|\ge2$, and its vertices
anchor both $(X_0)_1$ and $(X_0)_2$. Each has colour-zero
degree at most one into $P$, by unit-deficit insertion at $P_0$.

Let $k_c$ count the vertices of $M$ using colour $c=1,2$
towards $P$, and put $\gamma_c=a((X_0)_c)$, so that
$\gamma_1+\gamma_2=a_1+a_2\le p$. If $k_c\ge2$,
\cref{r3.12} makes all nonempty colour-$c$ neighbourhoods
in $P$ the same singleton. Its unit case also gives
$\gamma_c\ge2$. Thus $p\ge3$, and every vertex of $M$
must use the other colour $c'$, since its colour-zero and
colour-$c$ degrees are at most one. Applying the same argument
to $c'$ gives $\gamma_{c'}\ge2$ and degree at most one.
All three colour degrees are then at most one, so $p\le3$,
contrary to $p\ge\gamma_1+\gamma_2\ge4$.

Therefore $k_1,k_2\le1$. Since $p\ge a_1+a_2\ge2$, each
vertex of $M$ uses at least one of these two colours, giving
$|M|=2$ and $k_1=k_2=1$. The remaining two right vertices
must deviate from $(2,1)$ once in each row; the one deviating
in row one has vector $(*,1,1)$ and directly joins $P_1$ to
$(X_0)_1$. There are also $x\in P,y\in M$ with $xy$ of
colour one. These vertices give a second path
$P_1-x-y-(X_0)_1$, internally disjoint from the first, and
hence a two-core cycle.

After applying \cref{r3.2}, let the remaining types have sizes
$p\ge q\ge0$, and order the deficits
$\alpha_1\ge\alpha_2\ge\alpha_3$. If $\alpha_1\le p$,
the largest type is fully active and $p\le3$.
Put $d=5-R,x=p+w$. By (B), $Q$'s $d+1$ unit-column rows
have deficits at least $q+1$, so $dq+3\le x$; for $q=0$
this is $r_L\ge3$. The total budget gives $q\ge d+3-x$.
If $w\le d$, these force $p=3,w=d,q\le1$.
A unit row exists: the only alternative, $R=3$ with all
deficits two, is excluded by \cref{r3.29}. Then
\cref{r3.26} contradicts $p=3$. If $w>d$, separation is
full and (A) gives $p\le1$, hence $r_L=3$, impossible.
Therefore $\alpha_1\ge p+1$.

If $\alpha_2\le q$, (B) gives $R=4$ and two common active
rows. Their deficit sum is at most $q+w-1$. By (B), the nonunit
column contains one of these active cores of each type, so a
symmetry of the model takes $\boldsymbol b$ to $(2,1,1)$, the empty
type to $X_0$ and the common rows to $1,2$. Thus (C) forces
$w=2,p=q$ and full separation. A right permutation is
compatible only with the empty type; as in (A), all its edges to
$P\cup Q$ have its assigned colour $c$, so it is the centre of a
three-arm star in colour $c$. Under full separation the other member
of $Z_c$ is not a left permutation, and a vertex repeating $c$ would
close a two-core cycle with the star; so it is a right permutation
of the same type. By
\cref{r3.12}, $p=q=1$, say $P=\{u\},Q=\{v\}$.
Deleting $u,v$ separates the cores: a shortest remaining
core path avoids one star centre and closes with that star.
Right permutations occur in pairs and, by transposed
\cref{r3.23}, which applies as $r_L=p+q+w=4$, every $a_i\le2$ and
their type, having at least two vertices, is fully active after
transposition, number at most three, hence exactly two.
The exchange gives transposed full separation, contrary to (A).

Without right permutations, the four active cores force
$(1,2,1),(2,2,1)$ to occur, each at most once. Their forced
colour-zero edges enter active cores of unit right deficit,
so insertion gives $p=q=1$. These colour-zero cores have
no right anchors; $Z_0$ contains both repeated left vertices,
and the other separators are right trees. By \cref{r3.6},
the two cores require five left anchors, but only four exist.

Thus $\alpha_2\ge q+1$, forcing $R=w=3$ and deficits
$(p+1,q+1,1)$; \cref{r3.29} gives $p\ge2$.
We need one final observation: if a type of size at least two
has two active principal cores, the right vertices cannot all
anchor a unit-row core of its third colour $c$.
Indeed, (A) makes its third principal core inactive. Right
permutations can be assigned only $c$; the other separators
are repeated trees. By \cref{r3.11}, those colours have
degree zero when present in a right vector, and otherwise
degree at most one by unit-deficit insertion. By
\cref{r3.12}, two right vertices have no $c$ neighbour in
the type. They must miss both other colours, so are
constant-$c$ and give a two-core cycle.

Every column has an active core. If $q>0$, all right vertices
anchor the unit row's third core $D$, which is inactive.
Its column has an active core of type $V$. If $|V|\ge2$,
its unit-row and column-$D$ cores have the other two colours,
contrary to the observation. Otherwise $q=1,V=Q$, and
$P$ is active in the unit and deficit-two rows, violating (B).
If $q=0$, (A) makes the remaining $P_c$ inactive. Its column
has an active unit-row core; that row's last core has colour
$c$ and receives all right anchors, again contradicting the
observation. Undo the absorptions.
\end{proof}

\section{Permutation templates and finite transfer}\label{sec:permutations}

\subsection{A three-by-three Hall obstruction}\label{sec:6}
Let $I=C=\{0,1,2\}$ and identify $S_3$ with the bijections $I\to C$.
For $S\subseteq I\times C$, put
$w_c(S)=(|\{i:(i,c)\in S\}|-1)_+$ and $w(S)=\sum_cw_c(S)$.
A \emph{blocker} meets every permutation diagonal.

\begin{lemma}\label{r6.1}
Every blocker contains $H\times J$ with $|H|+|J|=4$.
It has weight zero if and only if it is one complete row.
Without a complete row its weight is at least two; without a complete row or column it contains a $2\times2$ rectangle.
\end{lemma}

\begin{proof}
Hall's theorem in the complement gives a nonempty $H$ with $|N(H)|<|H|$.
Choose $4-|H|$ columns outside $N(H)$.
The possible dimensions are $(1,3),(2,2),(3,1)$, and the latter two rectangles have weight two.
Weight zero allows at most one cell per column, so a complete row exhausts the blocker.
\end{proof}

\subsection{Capacity, absorption and realisation}\label{sec:7}
A \emph{core} is a nearly complete bipartite graph in a prescribed colour.
For each colour, the prescribed large-group edges form a union of vertex-disjoint cores.
Whenever an exceptional vertex has a colour label towards a group, it has at least $\alpha n$ neighbours of that colour in the group.
All core shores have size at least $\gamma n$, where $0<\alpha,\gamma\le1$ are fixed, and there are at most $m$ cores.
The homogeneous template completes these cores and makes each exceptional-to-group rectangle uniform in its label colour.

\begin{lemma}
\label{r7.1}\label{r7.2}\label{r7.3}
For sufficiently small $\delta>0$ and large $n$, the following hold when all prescribed missing degrees are at most $\delta n$.
\begin{enumerate}
\item A monochromatic cycle $C_0$ of half-length $h_0\le\delta n$ containing an edge of each selected core extends to a $C_{2n}$ if
\[
 \sum_DT_D\ge n-h_0,
 \qquad T_D=\min\{|A_D|-u_D,|B_D|-v_D\},
\]
where $u_D,v_D$ count its vertices in the two shores of $D$.
This inequality holds for two visited cores whose left shores have size at least $n-\delta n$ and whose right shores have total size at least $n$.
\item Let $L,R$ be prescribed anchors outside a core $A\times B$, each with at least $\alpha n$ neighbours of the core colour in the opposite shore, and let $|L|+|R|\le\delta n$.
All anchors can be absorbed into a $C_{2n}$ if $\min\{|A|+|L|,|B|+|R|\}\ge n$.
\item Suppose every pair of large groups is prescribed and every large group bearing a label colour belongs to a core of that colour.
A template $C_{2n}$ is realisable in the actual graph with the same exceptional vertices if there are at most $\delta n$ of them and their actual corner edges are retained.
\end{enumerate}
\end{lemma}

\begin{proof}
For (i), $T_D\ge\gamma n/2$.
Put $\beta=\min\{\gamma/4,1/(4m)\}$.
First assign $\lceil\beta n\rceil$ increments to every selected core, then distribute the remainder up to their capacities $T_D$.
Since $m\beta\le1/4$ and $4\delta<\beta$, this is possible and every positive increment satisfies \cref{r3.4}.
For the two-core case, the unused right capacities satisfy $x,y\ge\gamma n/2$ and $x+y\ge n-h_0$.
The smaller unused left capacity, truncated at $n$, is $N\ge n-2\delta n$.
Thus
\[
 \min(N,x)+\min(N,y)
 \ge\min(N+\gamma n/2,x+y)\ge n-h_0.
\]

For (ii), reserve two distinct neighbours per anchor and join the resulting paths cyclically inside the core by paths of length at most three.
The linear label degrees and missing-degree bound allow all $O(|L|+|R|+1)$ choices to avoid used vertices.
With no anchors, start with a four-cycle.
Its half-length $h_0$ leaves capacity
$\min\{|A|+|L|,|B|+|R|\}-h_0\ge n-h_0$.

For (iii), replace each maximal template core segment by an actual path of length at most three, with entries in the original label groups.
A singleton segment uses two entries in its shore and a common neighbour.
Choose all entries and internal vertices distinct, preserving the cyclic order and each segment's shore imbalance.
For $r$ exceptional vertices this uses $O_m(r+1)$ vertices, so greedy choices succeed for small $\delta$.
If $k_D^L,k_D^R$ and $u_D,v_D$ are the old and new core counts, then
\[
 \Delta_D=k_D^L-u_D=k_D^R-v_D\le T_D,
 \qquad \sum_D\Delta_D=n-h_0.
\]
A template cycle without exceptional vertices lies in one core with both shores of size at least $n$; start there with a four-cycle.
In (ii) and (iii), $h_0\le C_m(\delta n+1)$.
Choose $C_m\delta\ll\min\{\alpha,\gamma,\beta,1\}$.
Then $T_D\ge\gamma n/2$ and $m\lceil\beta n\rceil<n-h_0$ for large $n$, so the allocation above applies.
In (iii), $\sum_DT_D\ge\sum_D\Delta_D=n-h_0$ even if some $\Delta_D<0$.
\end{proof}

Since $\delta$ depends only on $m,\alpha,\gamma$, these assertions apply in particular when all missing degrees and exceptional counts are $o(n)$.

We return to a three-colouring of $K_{3n,3n}$.
All vertices outside the specified large groups are exceptional, so their counts equal the total deficits.

\begin{corollary}\label{r7.5}
Let $|A_i|=n-a_i$ and $|B_j|=n-b_j$, with all six deficits positive and their sum $o(n)$.
If colour $i+j$ has missing degree $o(n)$ on each $A_i\times B_j$, there is a monochromatic $C_{2n}$.
\end{corollary}

\begin{proof}
Take majority labels and retain the exceptional corner.
Since $n\ge3r+6$ for large $n$, \cref{r5.7} gives a template cycle, which \cref{r7.2}(iii) realises in the actual graph.
\end{proof}

\subsection{The forest on six permutation types}\label{sec:8}
\begin{theorem}\label{r8.2}
Fix $\gamma>0$.
Suppose $|A_i|=n-a_i$, where $a_i\ge1$ and $p=\sum_i a_i=o(n)$, and $|B_\pi|=\rho_\pi\ge\gamma n$ for every $\pi\in S_3$, with $\sum_\pi\rho_\pi=3n-q$ and $q=o(n)$.
Prescribe colour $\pi(i)$ on $A_i\times B_\pi$, with error degree $o(n)$.
If
\[
 \gamma n\le R_{ic}:=\sum_{\pi(i)=c}\rho_\pi\le(2-\gamma)n,
\]
then there is a monochromatic $C_{2n}$ for large $n$.
\end{theorem}

\begin{proof}
Suppose no such cycle exists.
Take majority labels; the homogeneous colour-$c$ cores are
$A_i\times D_{ic}$, where $D_{ic}=\bigcup_{\pi(i)=c}B_\pi$.
By \cref{r7.2}(iii), it suffices to find a template cycle.
Absorb exceptional left row vectors, namely $\pi\mapsto\pi(i)$, and right permutation vectors into their corresponding groups.
Recovery in \cref{r3.1} only uses $|A|+1>2r_{\rm opposite}$ and remains valid.
A left group reaching $n$ has a common colour neighbourhood of size at least $n$, giving a $K_{n,n}$.
Otherwise retain the notation; the two-sided bound on $R_{ic}$ persists with a smaller fixed $\gamma$.

Every colour incidence graph is a forest.
A cycle would lift to a short cycle through at least two cores, whose right shores have total size at least $(1+\gamma)n-q>n$, contrary to \cref{r7.1}(i).
For exceptional vectors $f_x:S_3\to C$ and $g_y:I\to C$, put
\[
 S_x=\{(i,c):\exists\pi,\ \pi(i)=f_x(\pi)=c\},
 \qquad S_y=\{(i,g_y(i)):i\in I\}.
\]
Let $L_{ic},R^+_{ic}$ count left and right incidences.
A cell is active if $L_{ic}\ge a_i$.
Absorption forces
\[
 d_{ic}:=n-R_{ic}\ge R^+_{ic}+1
\]
at every active cell: otherwise its left anchors and $\max(d_{ic},0)\le q$ right anchors fill both shores.
Since $\sum_cd_{ic}=\sum_cR^+_{ic}=q$, no row is wholly active.
Since $\sum_id_{ic}=q$, no column is wholly active when $q\le2$.
The forest charge bound in \cref{core:rules}(i) gives
\[
 W_c:=\sum_z w_c(S_z)\le2,
 \qquad W:=\sum_cW_c\le6.
\]
Each $S_x$ is a blocker other than a single row, and no $g_y$ is a permutation.
Thus every exceptional vertex contributes at least one to $W$, so $p+q\le6$.
If $q\ge3$, then $p=q=3$ and all $a_i=1$.
Every $S_x$ then lies in the active set and contains no complete row, giving $W\ge2p+q=9$, a contradiction.

Hence $q\le2$.
Each permutation diagonal has at least $p=\sum_i a_i$ left incidences, so the active set blocks $S_3$.
It has no complete row or column and contains $\{i,j\}\times\{c,d\}$ by \cref{r6.1}.
For the remaining row $k$ and colour $e$, the inactive cells $(k,c),(i,e),(j,e)$ meet at most $p-3$ left vertices.
Every other left vertex must meet both $(i,d)$ and $(j,d)$ to block the two diagonals through $(k,c)$.
Therefore
\[
 W_d\ge p-L_{kc}-L_{ie}-L_{je}\ge3,
\]
contrary to the bound $W_d\le2$.
Reverse the absorptions and realise the resulting cycle in the actual graph.
\end{proof}

\subsection{Arbitrary permutation support and finite transfer}\label{sec:9}
\begin{theorem}\label{r9.6}
Keep the group sizes and colour prescription of \cref{r8.2}, but allow any nonempty support $\mathcal F\subseteq S_3$ and drop the two-sided bound on $R_{ic}$.
The conclusion holds if $\rho_\pi\ge\gamma n$ for $\pi\in\mathcal F$, $\sum_{\pi\in\mathcal F}\rho_\pi=3n-q$, and $p,q,e=o(n)$, where $e$ is the maximum prescribed-colour error degree.
\end{theorem}

\begin{proof}
Suppose there is no monochromatic $C_{2n}$.
Use majority labels and cores $A_i\times D_{ic}$ as above, with left-anchor counts $L_{ic}$.
By \cref{r7.1}(ii), a cell with $L_{ic}\ge a_i$ absorbs once it has at least $\max(n-R_{ic},0)$ right anchors, since at most $p+q=o(n)$ anchors are involved.
A left-incidence cycle through cores of total right size at least $n$ also absorbs.
Indeed, choose distinct labelled entries and join them inside the cores.
The unused right capacities satisfy $x_D\ge\gamma n/2$, $\sum_Dx_D\ge n-h_0$, and the smallest unused left capacity, truncated at $n$, is $N=n-o(n)$.
For at least two visited cores,
\[
 \sum_D\min(N,x_D)
 \ge\min\{N+\gamma n/2,\sum_Dx_D\}\ge n-h_0,
\]
so \cref{r7.1}(i) applies.

If a row is wholly active, either some $R_{ic}\ge n$ and absorbs immediately, or the three positive right deficits sum to $q$ and the $q$ right majority labels fill one of them.
If $R_{i,\pi(i)}\ge n$ for all $i$, the $p$ left majority labels towards $B_\pi$ fill one left deficit.
In particular, every $\rho_\pi<n$ in a counterexample.

Call an opposite-parity pair strong when $\rho_\pi+\rho_\sigma\ge n$.
The strong-pair graph is a matching.
Otherwise some $\pi$ has strong neighbours $\sigma,\tau$, agreeing with it in distinct rows $j,k$; let $i$ be the third row.
These two cores have right size at least $n$, so their left-anchor counts are at most $a_j-1,a_k-1$.
At least $a_i+2$ left exceptional vertices avoid both; call them $G$.
All have label $\pi(i)$ towards $B_\pi$.
For $\nu\in\{\sigma,\tau\}$, let $r$ be its agreement row and $s$ the row other than $i,r$.
On $G$, the label towards $B_\nu$ is not $\pi(r)$.
At most one vertex has label $\pi(i)$: two would close a cycle through $(i,\pi(i))$ and $(s,\pi(i))$, since $\nu(s)=\pi(i)$, with right capacity at least $\rho_\pi+\rho_\nu\ge n$.
All others have label $\nu(i)=\pi(s)$.
Thus all three cells in row $i$ are active, a contradiction.

If a type $\delta$ is absent, give all other absent types weight zero.
Let $x_1,x_2$ be the remaining same-parity weights and $y_1,y_2,y_3$ the opposite weights.
Complete the strong matching in $K_{2,3}$ to one with edges $x_1y_2,x_2y_1$.
The other edges are weak, so
\[
 3n-q=x_1+x_2+y_1+y_2+y_3
 <2n+y_3<3n-\max(x_1,x_2).
\]
Hence $x_i<q<\gamma n$, forcing both to vanish.
Each $y_j<n$ and their sum exceeds $2n$, so all three survive.
They form one parity class, hence a Latin square after relabelling.
Their positive deficits sum to $q$, contradicting \cref{r7.5}.
Thus all six types survive.
Each has a weak opposite-parity neighbour of weight at least $\gamma n$, so $\rho_\pi<(1-\gamma)n$.
Consequently $2\gamma n\le R_{ic}<(2-2\gamma)n$, and \cref{r8.2} applies.
\end{proof}

For fixed $\gamma$, there are $\rho_P(\gamma),N_P(\gamma)>0$ such that $p,q,e\le\rho_P(\gamma)n$ and $n\ge N_P(\gamma)$ suffice in \cref{r9.6}.
Otherwise choose counterexamples with $n\to\infty$ and vanishing relative errors.

\begin{lemma}\label{finite:transfer}
For each $t>0$, choose $\rho(t)>0$ sufficiently small.
If $C_K\sqrt{s}<\rho(t)$, then for large $n$ the profile of \cref{r2.2} admits the following transfers simultaneously.
Call $V_x\times V_y$ an atomic rectangle; a union of atomic rectangles is uniquely coloured if each of its atomic pairs has a single allowed colour, and its edges of other colours are wrong.
Before cleaning, a uniquely coloured union rectangle has at most $C_Ksn^2$ wrong edges, and a rectangle containing an atomic pair on which colour $c$ is allowed has at least $tn^2$ edges of colour $c$.
Cleaning deletes, on each shore, the vertices of total wrong degree above $C_K\sqrt{s}n$; afterwards every allowed atomic colour retains at least $(t/2)n^2$ edges.
Three groups partitioning a shore, each of weight $h$, can be realised with sizes $n-a_i$, where $a_i>0$ and $\sum_i a_i\le C_K\sqrt{s}n$.
All prescribed error degrees are at most $C_K\sqrt{s}n$, and the anchors needed for all positive-colour extensions in \cref{r10.6} can be chosen together, pairwise distinct.
\end{lemma}

\begin{proof}
Each rectangle has at most $K^2$ atomic pairs.
Sum wrong edges over singleton-colour atomic pairs and delete, on each shore, vertices whose total wrong degree exceeds $C_K\sqrt{s}n$.
This loses at most $C_K\sqrt{s}n$ vertices and controls every prescribed rectangle.
The three group sizes were $n+O_K(s)n$; trim each to at most $n-1$, on both shores if necessary.
Each allowed atomic colour retains at least $(t/2)n^2$ edges, giving on either shore at least $(t/12)n$ vertices of degree at least $(t/12)n$.
There are at most $3K^2$ pools and nine target cores, with total anchor demand $O_K(\sqrt{s})n$.
Make this less than $(t/24)n$ and choose distinct anchors greedily.
Choose $\rho(t)$ also within the capacity tolerance for $m\le9$, $\alpha=t/12$, $\gamma=\min(t/4,1/4)$.
\end{proof}

Before regularity, fix $F$, small enough for \cref{r2.2}, so that
\[
 C_K\sqrt{F(t)}<\min\{t/100,\rho(t)/100,\rho_P(t/4)/100\},
 \qquad 3F(t)<t^2/100,
\]
and so that it meets any smaller strip tolerance below.
Call the resulting extracted profile \emph{admissible}.

\begin{lemma}\label{finite:strip}
An admissible profile has no partition $X=F_0\dotcup F_1\dotcup F_2$ and opposite set $D$ of weight greater than $h$ such that each $F_i\times D$ is entirely colour $i$.
The transposed assertion also holds.
\end{lemma}

\begin{proof}
The cover bound gives $w(F_i)=h$, and the extraction gap gives $w(D)\ge h+t/2$.
Its actual set has size at least $(1+t)n$.
Put omitted left vertices into the three rows, whose sizes remain $n+O_K(s)n$; the strip has $O_K(s)n^2$ wrong edges.
Delete vertices on either shore with total strip error degree above $C_K\sqrt{s}n$, then trim the rows to $n-a_i\le n$.
The common right set still has more than $n$ vertices, and exactly $\sum_i a_i$ left vertices lie outside the cores.
If some $a_i=0$, that core suffices.
Otherwise majority labels towards the common right set give at least $a_i$ anchors for some $i$, each of degree at least $n/3$.
Apply \cref{r7.3}(ii).
\end{proof}

A \emph{deterministic three-row permutation profile} is a partition $X=F_0\dotcup F_1\dotcup F_2$, with arbitrary weights, such that every $y\in Y$ has a permutation $\pi_y\in S_3$ for which $F_i\times\{y\}$ is entirely colour $\pi_y(i)$.

\begin{corollary}\label{r9.7}
An admissible profile is neither a deterministic three-row permutation profile nor its transpose.
\end{corollary}

\begin{proof}
A row heavier than $h$ forces its three colour neighbourhoods to have weight $h$, giving a transposed wide strip.
Thus all three rows have weight $h$.
By \cref{finite:transfer}, their positive deficits, exceptional counts and errors are $O_K(\sqrt{s})n$, while each present right type retains at least $(t/4)n$ vertices.
Apply the finite form of \cref{r9.6}, and transpose for the other case.
\end{proof}

We call the conclusion of \cref{r9.7} permutation exclusion and denote it by $(\mathrm P)$.

\section{Width and tripartition obstructions}\label{sec:obstructions}

\subsection{Strict width and degree rigidity}\label{sec:10}
Work in an admissible profile with shore weights $3h=1$, using the cover bound, \cref{finite:strip}, and $(\mathrm P)$; \cref{r10.6,r10.7} and \cref{r10.8} also use the transfer to the colouring given by \cref{r7.1,finite:transfer}.

\begin{lemma}\label{r10.1}
A two-colour support with shore weights $M,N\ge2h$ and component covers of weight at most $h$ is a deterministic checkerboard.
One shore has weight $2h$, split into two $h$-sets; the other is split into two sets, one possibly empty.
\end{lemma}

\begin{proof}
Both colours occur.
If one colour, say red, has one component, its cover $P\dotcup Q$ cannot meet both shores: otherwise blue contains $(X\setminus P)\times(Y\setminus Q)$ with both shore weights greater than $h$.
Transpose so that $Q=\varnothing$.
Blue contains $(X\setminus P)\times Y$, forcing $M=2h$, $w(P)=h$ and its cover to be $X\setminus P$.
Hence $P\times Y$ is red.

Otherwise neither colour has a zero-degree vertex, since its opposite-colour star would connect all components of the opposite colour.
Red has exactly two components; three would be connected by blue cross-rectangles.
The cross-rectangles are blue and the diagonals red, since a blue diagonal edge would join both blue blocks.
Each rectangle has a shore of weight at most $h$.
If a left group is heavier than $h$, both right groups have weight $h$; otherwise both left groups have weight $h$.
\end{proof}

\begin{proposition}\label{r10.2}
A component with one shore of weight $h$ has opposite shore weight at most $2h$.
In the equality case, its complementary two-colour rectangle is a checkerboard on four $h$-sets.
\end{proposition}

\begin{proof}
Let $A,B$ be red shores with $w(A)=h$ and $w(B)\ge2h$.
Apply \cref{r10.1} to $(X\setminus A)\times B$.
If $w(B)>2h$, the larger right checkerboard group $D$ is heavier than $h$.
Blue and green must be covered by their left $h$-groups, forcing $A\times D$ red, a forbidden strip.
The same contradiction excludes unequal right groups when $w(B)=2h$.
Unequal left groups force both right groups to have weight $h$; the larger left group is then red-joined to $Y\setminus B$, a transposed forbidden strip.
\end{proof}

For each partition into three $h$-sets on each shore, apply \cref{finite:transfer} once to obtain groups of sizes $n-a_i,n-b_j$.
Their total positive deficit and maximum error degree are at most $\eps_*n$, where $\eps_*=C_K\sqrt{s}$, small relative to the anchor margins.
Let $P,Q$ be the exceptional complements.
An anchor lies outside its core, possibly in another large group, and has linear degree into it.
A core is left- or right-active when anchors fill that deficit; activation on both sides gives a cycle.
A row is \emph{rainbow} if its three rectangles with the three columns are entirely coloured with three distinct colours, and similarly for a column.
An \emph{$h$-square} is an entirely coloured rectangle between two $h$-sets; a cover of weight at most $h$ of its component contains one of its two sides.

\begin{lemma}\label{r10.6}
The following rules hold simultaneously.
\begin{enumerate}
\item Entirely coloured rectangles give prescribed cores, and positive-colour extensions supply enough anchors to activate them.
\item A rainbow row has a right-active core, and a rainbow column has a left-active core.
\item Let rows $X_1,X_2$ be coloured $(B,G,R),(G,B,R)$ on columns $(Y_1,Y_2,T)$, and call the four cores $X_i\times Y_j$ with $i,j\in\{1,2\}$ checkerboard cores. Then some $Y_j$ has neither of its checkerboard cores left-active.
\end{enumerate}
\end{lemma}

\begin{proof}
Rule (i) is \cref{finite:transfer}.
For (ii), majority-label $Q$ towards the row.
The three counts sum to $|Q|=\sum_jb_j$, so one fills its deficit; transpose for a column.
For (iii), let $b_0,b_1,b_2$ be the deficits of $T,Y_1,Y_2$.
At most $b_0-1$ vertices of $Q$ have a red label towards either row, by absorption into $(X_1\cup X_2)\times T$.
For each other colour, at most one vertex has that label towards both rows: two give an eight-cycle through its two cores, extendible by \cref{r7.1}(i).
Every remaining vertex has distinct non-red labels and anchors both checkerboard cores in one column $Y_j$; assign it to $Q_j$.
If both columns had a left-active core, absorption would give $|Q_j|\le b_j-1$, and hence
\[
 |Q|\le(b_0-1)+(b_1-1)+(b_2-1)+2=|Q|-1.
\]
\end{proof}

\begin{lemma}\label{r10.7}
In a partition into three $h$-rows and three $h$-columns, if two rows are $(B,G,R)$ and $(G,B,R)$, the remaining corner is not red-free.
\end{lemma}

\begin{proof}
Name the rows $X_1,X_2,A$ and columns $Y_1,Y_2,T$, and suppose $A\times T$ is red-free.
Activate each checkerboard core with a non-red left extension from $A$.
Rule (iii) forces some $A\times Y_k$ entirely red.
Its rainbow column makes this red core left-active; a red extension to the other $Y_j$ would make it right-active, so $A\times Y_j$ is red-free.
A left extension of an $h$-square forces its cover to be the whole right shore.
Thus a non-red colour seen from $x\in A$ in $Y_j$ is forbidden from $x$ to $T$.
Both non-red colours cannot occur, since $A\times T$ is red-free.
Each $x\in A$ therefore has a deterministic permutation on $(Y_j,T,Y_k)$, as do $X_1,X_2$, contrary to transposed $(\mathrm P)$.
\end{proof}

\begin{theorem}\label{r10.8}
A component shore of weight $h$ has an opposite shore of weight less than $2h$.
\end{theorem}

\begin{proof}
Only equality in \cref{r10.2} remains.
Let $A,Y_1\dotcup Y_2$ be the red shores and $X_1,X_2$ the complementary checkerboard rows, all groups of weight $h$.
Put $T=Y\setminus(Y_1\cup Y_2)$.
Red is absent on $A\times T$ and positive on both $A\times Y_j$.
Let $L_{ij},R_{ij}$ denote non-red extensions of core $(i,j)$ into $A,T$.
An $h$-square can be covered only by a whole shore, so $L_{ij},R_{ij}$ cannot coexist, nor can same-colour squares join.
A blue left extension on $S\subseteq A$ forces $S\times T$ green, excluding every green right extension; the same holds with the colours exchanged.
By \cref{r10.7} and its transpose, extensions occur in both directions.
Hence all are blue and, after relabelling, only $L_{11},R_{22}$ occur.
Thus $A\times Y_2$ and $X_1\times T$ are red; $A\times Y_1$ and $X_2\times T$ use red and blue, with both on the former and blue on the latter.

Realise all activations together.
Blue $(1,1)$ is left-active, while blue $(2,2)$ and red $A\times Y_2$ are right-active.
If red occurs on $X_2\times T$, red $X_1\times T$ is left-active.
The rainbow column $Y_2$ and row $X_1$ then force green $(1,2)$ active on both sides.
Thus $X_2\times T$ is blue.
Its component has cover $X_2$, forcing $A\times T$ green.
Transpose with rows $(Y_2,T,Y_1)$ and columns $(A,X_1,X_2)$, and rename $(R,G,B)$ as $(B,G,R)$.
This contradicts \cref{r10.7}.
\end{proof}

We call the conclusion of \cref{r10.8} strict width and denote it by $(\mathrm B)$.

From here through the residual contradiction, the arguments are finite and use only the component-cover bound, wide-strip exclusion, permutation exclusion and strict width, together with \cref{r2.3,r2.4}, which need only rational weights and the cover bound.
Admissible profiles have these properties by \cref{r2.2,finite:strip,r9.7,r10.8}.

\begin{lemma}\label{r10.9}
If a colour has degree at least $h$ throughout one shore, its nonzero components are complete, their opposite shores have weight $h$, and all degrees on the given shore equal $h$.
Its critical components have only pure-sided minimum covers.
The assertion transposes.
\end{lemma}

\begin{proof}
Take a component with shores $U,V$, where $U$ lies in the given shore, and a cover $P\dotcup Q$ of weight at most $h$.
We have $w(V)\ge h$.
If $U\setminus P\ne\varnothing$, its degree bound gives $h\le w(Q)\le h-w(P)$, so $P=\varnothing$, $Q=V$ and $w(V)=h$.
Otherwise $w(U)\le h$.
Thus $w(V)>h$ implies $w(U)\le h$, while $w(V)=h$ implies $w(U)<2h$ by strict width.
If some $V$ were heavier than $h$, at most one other component would remain, and their left weights would total less than $h+2h=3h$.
Hence every $V$ has weight $h$, and the degree bound forces $W_c=1$ on $U\times V$.
Minimum covers of a complete rectangle are whole shores.
\end{proof}

\begin{lemma}\label{r10.10}
Two colours cannot both have degree at most $h$ throughout the same shore.
In particular, the three colour degrees cannot all be constantly $h$.
\end{lemma}

\begin{proof}
Transpose so that the given shore is $Y$.
The third colour has degree at least $h$, hence exactly $h$ by \cref{r10.9}; summing degrees forces all three equalities.
Their complete components have left shores of weight $h$ and right shores partitioning $Y$.
Each colour has two or three components by strict width; append its zero-degree left set if needed to obtain three $h$-groups.
Encode vertices by their three group indices in $\{0,1,2\}^3$.
Every opposite pair agrees in exactly one coordinate, and left coordinate marginals are uniform.
A right type has weight at most $h$ by strip exclusion, so at least three occur.

Let $a(u,v)$ count agreeing coordinates and normalise left weights $\lambda$ to sum to one.
For every $u$, $\sum_x\lambda_xa(u,x)=1$.
A present left type $u$ contributes $3\lambda_u>0$, forcing a present disjoint type $v$.
Relabel them $000,111$; right types are permutations of $012$.
Among three right types, two have the same parity and agree nowhere.
Relabel these $000,111$; left types are now permutations.
For every $\pi\in S_3$, uniform marginals give
\[
 1=3\lambda_\pi+
 \sum_{\operatorname{sgn}(\sigma)\ne\operatorname{sgn}(\pi)}\lambda_\sigma.
\]
Weights are therefore equal within parity classes, including zero weights for absent types.
If only one class is positive, its three left types relabel as $000,111,222$, giving $(\mathrm P)$.
If both are positive, all six left permutations occur.
A nonconstant right type agrees twice or three times with some left permutation, so right types must be the three constants, giving transposed $(\mathrm P)$.
\end{proof}

\subsection{Critical-side tripartitions}\label{sec:11}

\subsubsection{Structural assumptions}
Continue with a finite profile on shores of weight $3h$ satisfying the component-cover bound, wide-strip exclusion $(\mathrm S)$, permutation exclusion $(\mathrm P)$ and strict width $(\mathrm B)$.
Here $(\mathrm S)$ asserts that there is no partition $X=F_0\dotcup F_1\dotcup F_2$ and right set $D$ with $w(D)>h$ such that every $F_i\times D$ is entirely colour $i$, and the same holds after interchanging the shores.
Selected critical components need not initially have only pure-sided covers.

Suppose each colour $i$ has a selected critical component with shores $A_i,B_i$, where
\[
 X=A_0\dotcup A_1\dotcup A_2,\qquad w(A_i)=h.
\]
By \cref{r2.12},
$B_i=D\dotcup E_i$, where $D=B_0\cap B_1\cap B_2$, the four sets $D,E_0,E_1,E_2$ are pairwise disjoint, and $A_i\times D$ is entirely colour $i$.
Put
\[
 d=w(D),\qquad e_i=w(E_i),
 \qquad T=Y\setminus(D\cup E_0\cup E_1\cup E_2).
\]
Strip exclusion and criticality give $d\le h$ and $d+e_i\ge h$.
We do not yet assume purity or apply $(\mathrm B)$ to the selected components.
For distinct $i,j,k$, the rectangle $A_j\times E_k$ has colour $i$, since the other colours are excluded by the selected supports.

\subsubsection{The common part}
\begin{lemma}\label{r11.2}
The three private-part weights $e_0,e_1,e_2$ cannot all be positive.
\end{lemma}

\begin{proof}
Suppose all three are positive.
Call colour $i$ joining if $A_j\times E_k$ and $A_k\times E_j$ lie in one component, and separating otherwise.
A separating colour has whole left sides $A_i,A_j,A_k$, hence right degree at most $h$ and neighbours in at most one row at each point of $T$.
By \cref{r10.10}, the joining set $J$ has size $r\ge2$.
A joining component's cover must contain $E_j\cup E_k$: omitting part of either forces one whole row into the cover, leaving the other positive rectangle uncovered.
Write the rest as $S_i\dotcup T_i$, where
\[
 S_i\subseteq A_j\cup A_k,\qquad T_i\subseteq T,
 \qquad w(S_i)+w(T_i)\le\lambda_i=h-e_j-e_k.
\]
For $i\in J$, since $A_j\times E_k$ and $A_k\times E_j$ are complete and nonempty, the joining colour-$i$ component has left shore exactly $A_j\cup A_k$.
A point of $T$ lies outside $B_i$, so all its colour-$i$ neighbours lie in this component, and those of a point outside $T_i$ lie in $S_i$.
Put $s=\sum_ie_i$ and $t=w(T)=3h-d-s$.
Outside $\bigcup_{i\in J}T_i$, separating colours meet at most $3-r$ rows, so $r$ whole rows must lie in the $S_i$.
This gives $rh\le\sum_{i\in J}w(S_i)<rh$.
Thus $T=\bigcup_{i\in J}T_i$ and $t\le\sum_{i\in J}\lambda_i$.
Since $t-\lambda_i-\lambda_j=h-d+e_k>0$, all three colours join.
Their budgets give $\sum_iw(S_i)\le d-s<h$.
A point in only one $T_i$ would require the other two $S$-sets to cover $A_i$.
Hence
\[
 2t\le\sum_iw(T_i)\le3h-2s,
 \qquad 2t-(3h-2s)=3h-2d>0,
\]
a contradiction.
\end{proof}

\begin{proposition}\label{r11.3}
We have $d=h$, and the three selected components have no mixed minimum vertex cover.
\end{proposition}

\begin{proof}
If $d<h$, then $d+e_i\ge h$ makes all $e_i$ positive, contrary to \cref{r11.2}.
Each $A_i\times D$ is therefore a complete rectangle with both shore weights $h$.
Any cover contains all of $A_i$ or all of $D$, so a minimum cover of the selected component is one of these pure-sided sets.
Now $(\mathrm B)$ gives $0\le e_i<h$.
\end{proof}

\subsubsection{Cover capacity}
\begin{proposition}\label{r11.7}
A critical-side tripartition with common right part of weight $h$ and $0\le e_i<h$ cannot exist.
\end{proposition}

\begin{proof}
Call an additional colour-$i$ component left-wide if its left shore $L_i\subseteq X\setminus A_i$ has weight greater than $h$.
There is at most one per colour.
Let $J$ index them and put $j=|J|$.
Colours outside $J$ have right degrees at most $h$, so \cref{r10.10} gives $j\ge2$.
For $i\in J$, choose arbitrary covers $P_i\dotcup Q_i$ with $p_i+q_i\le h$, and put $R_i=X\setminus(A_i\cup L_i)$.
Then $p_i<h$, $q_i>0$ and $w(R_i)<h$.
For $i\notin J$, set $P_i=Q_i=L_i=\varnothing$.

For $i\in J$, a right vertex outside $B_i\cup Q_i$ has its colour-$i$ neighbourhood contained in $P_i$ or $R_i$, hence degree less than $h$.
On $B_i$ this degree is at most $h$; colours outside $J$ have degree at most $h$ everywhere.
The three degrees sum to $3h$, and any two $B_i$ meet in $D$, so
\[
 \bigcup_iQ_i=Y\setminus D,
 \qquad 2h\le\sum_iq_i\le jh-p,
 \qquad p:=\sum_ip_i\le h.
\]
Put $U_i=Q_i\setminus\bigcup_{k\ne i}Q_k$, $I=\{i:w(U_i)>0\}$ and $u=\sum_iw(U_i)$.
Whenever $k\in J$ and $k\ne i$, the private rectangle gives colour-$k$ degree at least $h$ on $E_i$.
Since $E_i\cap B_k=\varnothing$, the preceding degree bound gives $E_i\subseteq Q_k$.

We claim $|I|+j\le4$.
Otherwise $j=3$ and two $U_i$ have positive weight; the preceding inclusions give $E_i\subseteq Q_j\cap Q_k$ and hence $U_i\subseteq T$.
At a point of $U_i$, the neighbourhoods cover $X$ using $L_i$ and one of $P_k,R_k$ for each other colour.
For distinct such $i,j$, write $P=\bigcup_kP_k$ and let $k$ be the remaining index.
Then
\[
 X\setminus(P\cup R_k)
 \subseteq(L_i\cup R_j)\cap(R_i\cup L_j)
 =(L_i\cap L_j)\cup(R_i\cap R_j).
\]
The last set avoids $A_i\cup A_j$, forcing $2h\le p+w(R_k)<2h$.

Points outside the unique parts of the $Q_i$ have multiplicity at least two, so $\sum_iq_i\ge4h-u$.
Thus
\[
 (4-j)h+p\le u\le\sum_{i\in I}q_i
 \le |I|h\le(4-j)h.
\]
Hence $p=0$ for every choice of covers of weight at most $h$, and each $Q_i$, $i\in J$, is the whole right shore of its component.
Equality also gives $|I|=4-j$, and $Q_i=U_i$, $q_i=h$ for $i\in I$.
If $j=2$, the two $Q_i$ partition $Y\setminus D$ into $h$-sets.
If $j=3$, relabel so that $I=\{0\}$.
Then $Q_0$ is disjoint from $Q_1\cup Q_2$, and $U_1=U_2=\varnothing$ gives $Q_1=Q_2=(Y\setminus D)\setminus Q_0$.
Thus all $q_i=h$ for $i\in J$; the left-wide components are pure-sided critical components with left weights in $(h,2h)$ by strict width.

Write $F_0=Q_0$ and $F_1=(Y\setminus D)\setminus F_0$, relabelling so that $Q_1=F_1$; every other nonempty $Q_i$ equals $F_1$.
Applying $E_i\subseteq Q_k$ only for $k\in J$, $k\ne i$, and using $Q_0\cap Q_1=\varnothing$, gives
\[
 E_2=\varnothing,\qquad E_0\subseteq Q_1,
 \qquad E_1\subseteq Q_0.
\]
In particular, $B_2=D$.
Apply \cref{r11.3} after transposing to the critical right-side partition $Q_0,Q_1,D=B_2$.
Its common left part is $A_2$; the common-rectangle argument and strict width give
\[
 L_0=A_2\dotcup V_1,\qquad L_1=A_2\dotcup V_0,
 \qquad V_i\subseteq A_i,\qquad0<w(V_i)<h.
\]
Colour $2$ contains
\[
 V_0\times Q_0,\quad (A_1\setminus V_1)\times(Q_0\setminus E_1),
 \quad V_1\times Q_1,\quad (A_0\setminus V_0)\times(Q_1\setminus E_0).
\]
The first pair lies in one component and the second in another, since $e_i<h$ makes the shared right sets nonempty.
If these components coincided, any cover $P\dotcup Q$ of weight at most $h$ would contain $V_0,V_1$ in $P$; omitting either forces its whole column into $Q$ and leaves the other positive rectangle uncovered.
It cannot also contain a whole complementary $A_i\setminus V_i$, since that would exceed $h$.
The displayed rectangles and $A_0\times E_1,A_1\times E_0$ therefore force both
$(Q_1\setminus E_0)\cup E_1$ and $(Q_0\setminus E_1)\cup E_0$ into $Q$, whose union has weight $2h$.
Thus the components are distinct.
If $E_0$ or $E_1$ were nonempty, its complete join to $A_1$ or $A_0$ would connect them.
All private parts therefore vanish, and the colour-$2$ components have whole right shores $Q_0,Q_1$.

Every colour now has, besides its component on $D$, a component with whole right shore $F_0$ or $F_1$.
For $x\in A_i$, colour $i$ occupies $D$; each other colour is restricted to one of $F_0,F_1$ by its additional component.
They must fill distinct whole columns, giving the transpose of a deterministic three-row permutation profile, contrary to $(\mathrm P)$.
\end{proof}

\begin{theorem}\label{r11.8}
Under the component-cover bound, $(\mathrm P)$, $(\mathrm B)$ and $(\mathrm S)$, one cannot choose a critical component in each colour whose whole shores of weight $h$ partition the same shore.
The selected components may initially have mixed minimum vertex covers.
\end{theorem}

\begin{proof}
By \cref{r11.3}, the common opposite part has weight $h$ and the selected components have only pure-sided minimum covers.
Strict width gives $0\le e_i<h$, contradicting \cref{r11.7}.
The transposed case is identical.
\end{proof}

\section{Critical components and the cover partition}\label{sec:critical}

\subsection{Two critical components in each colour}\label{sec:12}

Both shores have weight $3h$. We use the balanced covers in
\cref{r2.4}, strict width $(\mathrm B)$, and the obstructions in
\cref{r10.10,r11.8}. By \cref{r2.3}, each colour has between one
and three critical components. A minimum cover contained in one
shore is the whole corresponding component side, of weight $h$;
this does not yet imply purity or uniqueness.

\begin{lemma}\label{r12.1}
Let $E_1\times F_1,E_2\times F_2$ be complete rectangles in one
colour, where $E_1,E_2$ are disjoint sets of weight $h$ and
$w(F_1),w(F_2)>0$. If $w(F_1\cup F_2)\ge h$, they belong to
different components. Consequently, if three left $h$-sets partition
$X$, each is completely joined in this colour to a nonempty right
set, and every pair of right sets has union of weight at least $h$,
then every right degree in this colour is at most $h$.
\end{lemma}
\begin{proof}
Suppose that the rectangles belong to one component, with a cover
$H$ of weight at most $h$. Omitting a vertex of $F_i$ forces
$E_i\subseteq H$, using the whole budget, while the other rectangle
still needs positive cover weight. Thus $F_1\cup F_2\subseteq H$,
so $w(F_1\cup F_2)=h$ and $H=F_1\cup F_2$. This is the whole
right side, while the left side has weight at least $2h$, contrary
to $(\mathrm B)$. For the final assertion, the three rectangles
lie in distinct components. Their left sides are exactly the three
specified sets, since these partition $X$, giving the degree bound.
\end{proof}

\begin{lemma}\label{r12.2}
Let $X=E_0\dotcup E_1\dotcup E_2$, with $w(E_i)=h$. Suppose
$E_0,E_1$ are whole left sides of critical components in different
colours. If the third colour contains $E_2\times F$ completely,
where $w(F)\ge h$, then $w(F)=h$. Its component has whole right
side $F$ and left weight greater than $h$.
\end{lemma}
\begin{proof}
The rectangle requires cover weight at least $h$, so its component
is critical. Its whole left side cannot be $E_2$, by \cref{r11.8}.
A cover of weight at most $h$ therefore cannot contain $E_2$,
since the additional left vertices still need their edges covered.
It must contain $F$. Hence $w(F)=h$, the unique minimum cover is
$F$, and $F$ is the whole right side. The left side properly
contains $E_2$.
\end{proof}

\begin{lemma}\label{r12.3}
Three critical components involving at least two colours cannot
have whole sides of weight $h$ that partition one shore. No
restriction is imposed on other critical components or on the
orientations or purity of the selected components' minimum covers.
\end{lemma}
\begin{proof}
The three-colour case follows from \cref{r11.8}. Otherwise,
interchange the shores if necessary, rescale to $h=1$, and suppose
that a red side $A$ and blue sides
$B_1,B_2$ partition $X$. Their right sides $S,S_1,S_2$ have
weights in $[1,2)$, with $S_1\cap S_2=\varnothing$. Every blue
component has left weight at most one, since all remaining blue
left sides lie in $A$. Thus all blue right degrees are at most one.
Green contains $B_1\times(S\setminus S_1)$ and
$B_2\times(S\setminus S_2)$. Their right sets have union $S$;
if they intersected, they would belong to one component, contrary
to \cref{r12.1}. Thus $S\subseteq S_1\cup S_2$. Put
$D_i=S\cap S_i$, $d_i=w(D_i)$, and
$U=(S_1\cup S_2)\setminus S$, $u=w(U)>0$.
By \cref{r12.2}, $d_i\le1$.

The green rectangles $A\times U,B_1\times D_2,B_2\times D_1$
have disjoint right sets, with
\[
 d_1+d_2\ge1,\qquad
 u+d_i=w(S_1)+w(S_2)-d_{3-i}\ge1.
\]
If all three sets were nonempty, \cref{r12.1} would bound green
right degrees by one, contradicting \cref{r10.10}. Since $U$ is
nonempty and $w(D_1\cup D_2)\ge1$, exactly one $D_i$ is empty.
Relabel so that $D_1=\varnothing$ and $D_2=S$.
Then $w(S)=1$ and $u\ge1$.
The green components $G_0,G_1$ containing $A\times U,B_1\times S$
are distinct. By \cref{r12.2}, $G_1$ has whole right side $S$
and left weight greater than one.

The whole left side of $G_0$ is $A$. Otherwise its cover must be
$U$, forcing $u=1$ and $w(S_1)=w(S_2)=1$. Blue then has whole
right sides $S_1,S_2$ of weight one, and green has $U,S$; every
other right side lies in the corresponding weight-one complement.
Both colours would have all left degrees at most one, contrary
to \cref{r10.10}.

Thus $B_2\times S_1$ is red: blue is forbidden by component
separation, and green by $S_1\subseteq U$ and the left side $A$
of $G_0$. Apply \cref{r12.2} to the green side $A$, blue side
$B_1$, and this rectangle. Its red component has whole right
side $S_1$ and left weight greater than one. Its left side avoids
$A$, so it meets $B_1$ in a positive-weight set $V$. Red and blue
neighbours of $V$ lie in $S_1$, and green neighbours lie in $S$.
Since $w(S\cup S_1)=2<3$, this is impossible.
\end{proof}

\begin{lemma}\label{r12.4}
Let $w(X)=3h$ and $w(A)=h$. Let each of the families $\{B_i\}$
and $\{G_j\}$ consist of at most two disjoint sets of weight $h$. If
nonnegative coefficients, with $\alpha>0$, satisfy
\[
 \alpha\ind{A}+\sum_i\beta_i\ind{B_i}
 +\sum_j\gamma_j\ind{G_j}=\ind{X},
\]
then $A$ and two sets from these families partition $X$. The two
additional sets may belong to the same family.
\end{lemma}
\begin{proof}
Discard zero coefficients and put $Z_B=X\setminus\bigcup_iB_i$,
so $w(Z_B)\ge h$. If $Z_B\subseteq A$, then $Z_B=A$ and two
$B$-sets partition $X\setminus A$. Otherwise, at a vertex of
$Z_B\setminus A$ only one $G$-term remains, forcing some
$\gamma_j=1$. Symmetrically, either two $G$-sets give the required
partition with $A$, or some $\beta_i=1$. These coefficient-one
sets are disjoint from each other and from $A$, since $\alpha>0$.
Their weights are $h$, so the three sets partition $X$.
\end{proof}

\begin{lemma}\label{r12.5}
For each colour choose at most two whole left sides $A_{cj}$ of
critical components, each of weight $h$. There is no identity
\[
 \sum_{c,j}\alpha_{cj}\ind{A_{cj}}=\ind{X},
 \qquad\alpha_{cj}\ge0,\qquad
 \sum_j\alpha_{cj}>0\quad(c=0,1,2).
\]
The transposed assertion also holds.
\end{lemma}
\begin{proof}
Discard zero coefficients. Two colours cannot each retain two
sides: the two selected sides use weight $2h$ and leave only $h$.
Every component of either colour would have left weight at most
$h$, giving right degrees at most $h$ in both colours,
contrary to \cref{r10.10}. Thus one colour retains one side $A$
with positive coefficient. By \cref{r12.4}, $A$ belongs to a
partition into three critical sides involving at least two
colours, contrary to \cref{r12.3}. Interchange the shores for the
transposed assertion.
\end{proof}

\begin{lemma}\label{r12.6}
Choose a minimum cover of every critical component. If each colour
has at most two chosen covers entirely on the left and at most
two entirely on the right, their union is $X\dotcup Y$.
\end{lemma}
\begin{proof}
Use \cref{r2.4}, whose coefficient sum is positive in every colour.
Taking weights gives
$h\sum_{c,j}a_{cj}=3h(k_X+k_Y)>0$. If $k_X=0$, all covers with
positive coefficient are entirely right, and $k_Y>0$. Normalising
by $k_Y$ gives the identity excluded by the transpose of
\cref{r12.5}. Similarly, $k_Y\ne0$. Thus both constants are
positive, and every vertex belongs to a chosen cover.
\end{proof}

\begin{theorem}\label{r12.7}\label{r12.8}
Each colour has exactly two critical components. Every choice of
their six minimum covers partitions $X\dotcup Y$. Consequently,
every critical component has a unique minimum cover.
\end{theorem}
\begin{proof}
Suppose red has three critical components. Their left sides
$A_1,A_2,A_3$ and right sides $D_1,D_2,D_3$ partition the shores
into $h$-sets. Another such colour would give two colours with
right degrees at most $h$, contrary to \cref{r10.10}. Thus blue
and green have at most two critical components each. Fix their
minimum covers, with union $E_X\dotcup E_Y$. Then
\[
 w(E_X)+w(E_Y)\le4h.
\]
For each $i$, choose red covers $A_i$ and $D_j$ for $j\ne i$.
By \cref{r12.6},
\[
 X\setminus A_i\subseteq E_X,\qquad D_i\subseteq E_Y.
\]
The blue and green covers stay fixed. Taking all three $i$ forces
$E_X=X,E_Y=Y$, contradicting $w(E_X)+w(E_Y)\le4h$.

Every colour therefore has at most two critical components.
Arbitrary minimum covers satisfy \cref{r12.6}: at most six
$h$-sets cover weight $6h$. There are exactly six and they
partition the space. Fixing five covers makes the sixth the
complement of their union. Since the choice was arbitrary, this
proves uniqueness.
\end{proof}

\subsection{The six-cover partition}\label{sec:13}
Rescale to $h=1$, so each shore has weight three.

\begin{proposition}\label{r13.1}\label{r13.2}\label{r13.3}
At least two critical covers are mixed, and at most two pure
covers lie on each shore. A critical component side has weight
one exactly when it is the minimum cover. The side opposite a
pure cover has weight in $(1,2)$, while both sides of a component
with a mixed cover have weight greater than one. For the residual
sets $T_c^X,T_c^Y$ outside the two critical components of colour
$c$, distinct colours satisfy
\[
 w(T_c^X)+w(T_d^X)<2,\qquad w(T_c^Y)+w(T_d^Y)<2.
\]
\end{proposition}
\begin{proof}
Use the partition and uniqueness in \cref{r12.8}. A unit side is
a minimum cover, proving the criterion; two unit sides in one
component contradict uniqueness. Strict width gives the opposite
upper bound. If all covers were pure, three critical sides would
partition a shore, contrary to \cref{r12.3}. A single mixed cover
cannot give integer total left weight three: its nonintegral left
weight is added to integer contributions from the pure covers.
Three pure covers
on one shore would exhaust it, excluding mixed covers.

Each residual weight is at most one. Equality for two colours
would give two unit critical left sides in each colour. All their
components would then have left weight at most one, contradicting
\cref{r10.10}. Interchange the shores for the right-residual case.
\end{proof}

\section{Residual sets and the main theorem}\label{sec:14}\label{sec:residual}

\subsection{Notation}
Use \cref{r12.8,r13.3}. Write the $j$th critical cover of colour
$c$ as
\[
 H_{cj}=P_{cj}\dotcup Q_{cj},\qquad
 P_{cj}\subseteq X,\quad Q_{cj}\subseteq Y,\qquad
 p_{cj}+q_{cj}=1,
\]
where lower-case letters denote weights. The two critical covers
of one colour are \emph{companions}. Let $T_c^X,T_c^Y$ be
the sets outside the two critical component sides of colour $c$,
and put $a_c=w(T_c^X)$, $b_c=w(T_c^Y)$. Both are at most one.
Residual sets may be disconnected; only component supports will
be used. By \cref{r13.1},
\[
 a_c+a_d<2,\qquad b_c+b_d<2\qquad(c\ne d);
\]
we refer to these inequalities as the residual budget.

\subsection{A pure pair forced by crossing covers}

\begin{lemma}\label{r14.1}
A crossing rectangle $P_{c1}\times Q_{c2}$ or
$P_{c2}\times Q_{c1}$ uses only the other colours' residual
components. If a nonempty rectangle $E\times F$ uses only blue
and green residual components, then
\[
 E\subseteq T_B^X\cup T_G^X,\qquad
 F\subseteq T_B^Y\cup T_G^Y.
\]
If $E\not\subseteq T_B^X$, then $F\subseteq T_G^Y$; the
symmetric assertions also hold.
\end{lemma}
\begin{proof}
Colour $c$ cannot join its distinct critical components. The
endpoints lie in no other chosen cover, so another colour cannot
use a critical component either. The residual support assertions
follow because every opposite pair has an allowed colour, giving
the displayed inclusions. In particular, a vertex of $E\setminus T_B^X$ can meet $F$ only
through green residual components.
\end{proof}

\begin{lemma}
\label{r14.pure}\label{residual:forcing}
\begin{enumerate}
\item A pure left cover and a pure right cover exist.
\item Some colour has two pure covers on the same shore.
\end{enumerate}
\end{lemma}
\begin{proof}
For (i), the residual budget gives $\sum_ca_c<3$. If there were no
pure left cover, each nonempty $P_{ci}$ would have a companion
with nonempty right part $Q_{c,3-i}$. By \cref{r14.1}, all left-cover parts,
and hence $X$, would lie in the three left residual sets, a
contradiction. Interchange the shores for a pure right cover.

For (ii), suppose no colour has two pure covers on one shore.
All residual weights are then less than one: equality would leave
two unit component sides, which are pure covers by \cref{r13.3}.
Choose a pure red left cover $A$ using (i), with companion
$P_R\dotcup Q_R$. We have $q_R>0$, since a second pure red left
cover is excluded. As $w(A)=1>a_B,a_G$, \cref{r14.1} gives
\[
 Q_R\subseteq T_B^Y\cap T_G^Y.
\]
Thus $q_R<1$ and $0<p_R<1$. Another pure left cover, say blue,
would similarly have a mixed companion $P_B\dotcup Q_B$ with
$Q_B\subseteq T_G^Y$. But $p_R+p_B\le1$, while these disjoint
right parts in $T_G^Y$ have weight $2-p_R-p_B\ge1>b_G$.
Hence $A$ is the only pure left cover.

For $c\in\{B,G\}$, let $d$ be the other colour. Every nonempty
$P_{ci}$ has a nonempty companion right part, so
$P_{ci}\subseteq T_R^X\cup T_d^X$ by \cref{r14.1}. If
$x\in P_{ci}\setminus T_R^X$, then $x\in T_d^X$. It lies
outside $A$ and all red covers, so it is in the mixed red
component outside $P_R$. Its red neighbours lie in
$Q_R\subseteq T_d^Y$, as do its colour-$d$ neighbours. Thus
$\{x\}\times(Y\setminus T_d^Y)$ is entirely colour $c$,
putting weight $3-b_d>2$ in its critical right side, while the
other critical right side has weight at least one. This is
impossible. All blue and green left-cover parts therefore lie
in $T_R^X$, forcing
$a_R\ge w(X\setminus(A\cup P_R))=2-p_R>1$.
\end{proof}

\begin{lemma}\label{r14.2}
Let $P\dotcup Q$ be a mixed minimum cover of a red critical
component. The inclusions
\[
 P\subseteq T_B^X\cap T_G^X,\qquad Q\subseteq T_B^Y\cap T_G^Y
\]
cannot both hold.
\end{lemma}
\begin{proof}
Suppose both hold, and put $p=w(P)$, $q=w(Q)$, so $p,q>0$ and
$p+q=1$. The rectangle
\[
 \bigl(X\setminus(T_B^X\cup T_G^X)\bigr)\times Q
\]
is entirely red. Its left set is nonempty by the residual budget and
belongs to the red component containing $Q$ and $P$. This
component has left weight at most two, since the other red
critical left side has weight at least one. The two specified
left subsets are disjoint, and $w(T_B^X\cap T_G^X)\ge p$. Hence
\[
 2\ge p+3-w(T_B^X\cup T_G^X)\ge3-a_B-a_G+2p.
\]
Similarly, $b_B+b_G\ge1+2q$. With the residual budget, this gives
\[
 4=2+2(p+q)\le a_B+a_G+b_B+b_G<4,
\]
a contradiction.
\end{proof}

\subsection{Excluding the pure pair}

\begin{lemma}\label{r14.3}\label{r14.4}
The two critical covers of a colour cannot both be pure on the
same shore. Consequently, all six residual sets have weight less
than one.
\end{lemma}
\begin{proof}
By symmetry, suppose both green covers are pure left, with whole left sides
$A_1,A_2$. Put $A_0=X\setminus(A_1\cup A_2)=T_G^X$ and $Z=T_G^Y$.
Then $w(A_0)=1$ and $w(Z)<1$. Neither other colour has a pure
left cover, so their four critical left sides have weight greater
than one. By \cref{r14.pure}(i), choose a pure red right cover
$D_R$, and write the other red cover as $P_R\dotcup Q_R$.
All red and blue left-cover parts lie in $A_0$. The crossing
$P_R\times D_R$ gives $P_R\subseteq T_B^X$, since
$w(D_R)=1>w(Z)$.

Write the blue covers as $P_i\dotcup Q_i$. We claim
$P_i\subseteq T_R^X$ for both $i$. Otherwise, with $j\ne i$,
\cref{r14.1} gives $Q_j\subseteq Z$. In fact,
$Q_j\subseteq T_R^Y$: at $y\in Q_j\setminus T_R^Y$, green
neighbours lie in $A_0$, and red neighbours lie in the red
left-cover part, also in $A_0$. Thus $A_1\cup A_2$ is completely
blue-joined to $y$, forcing one blue critical left side to have
weight at least two, while the other has weight greater than one.
Now
\[
 w(Q_1)+w(Q_2)=1+w(P_R)>w(Z),
\]
so $Q_i\not\subseteq Z$. Here $Q_j$ is nonempty because blue
has no pure left cover. \cref{r14.1} gives
$P_j\subseteq T_R^X$. Since $0<w(Q_j)<1$, this cover is mixed
and lies in both red and green residual sets, contrary to
\cref{r14.2}. This proves the claim.

Both blue left-cover parts lie in the red and green residual sets.
If both were nonempty, their complete blue
joins to $Y\setminus(T_R^Y\cup Z)$ would connect their critical
components; this common right set has weight at least
$3-b_R-b_G>1$ by the residual budget. Their total weight is
$1-w(P_R)>0$, and neither has weight one, since no red or blue
cover is pure left. Thus exactly one blue left-cover part is
nonempty, say $P_B$, with $0<w(P_B)<1$. Blue has a mixed cover
$P_B\dotcup Q_B$ and a pure right cover $D_B$, and
\[
 A_0=P_R\dotcup P_B,\qquad
 Y=D_R\dotcup D_B\dotcup E,\qquad
 E=Q_R\dotcup Q_B,\quad w(E)=1.
\]
The red cover $P_R\dotcup Q_R$ is mixed as well.

Let $L_R,L_B$ be the whole left sides corresponding to $D_R,D_B$.
Both have weight greater than one. They avoid their own mixed
components, and $P_R\subseteq T_B^X$, $P_B\subseteq T_R^X$,
so $L_R,L_B\subseteq A_1\cup A_2$. Thus
$I=L_R\cap L_B$ has positive weight, and $I\times E$ is green.
If $I$ met both $A_1,A_2$, it would join the green components.
After relabelling, $I\subseteq A_1$ and $E$ lies in that green
component. The set $V=L_R\cap A_2$ has positive weight.
On $V\times E$, red is forbidden by the component with cover
$D_R$, and green by component separation. Hence this rectangle
is blue. Since $Q_B$ is nonempty, $V$ belongs to the mixed blue
component. But $V\cap P_B=\varnothing$, so its blue neighbours
must lie in $Q_B$, contradicting the blue rectangle $V\times Q_R$.

Finally, a residual side of weight one leaves two unit critical
component sides. Both would be pure covers on that shore, which
we have excluded.
\end{proof}

\begin{theorem}\label{r14.5}
Under the component-cover bound, wide-strip and permutation
exclusions, and strict width, the critical-component count
cannot be $(2,2,2)$.
\end{theorem}
\begin{proof}
The six-cover structure follows from \cref{r12.8,r13.3}.
The pure pair given by \cref{residual:forcing}(ii) is excluded
by \cref{r14.3}.
\end{proof}

\subsection{Proof of the finite structural alternative}\label{sec:structural-proof}
\begin{proof}[Proof of \cref{struct:alternative}]
Suppose none of the three configurations occurs. The third
exclusion gives strict width. In any wide strip, the cover bound
forces the three row weights to be at most $h$, hence all equal
to $h$. An unequal-weight three-row permutation profile has a
row heavier than $h$; its three colour neighbourhoods partition
the opposite shore and each have weight at most $h$, hence
exactly $h$. This gives a transposed wide strip. Thus the first
two exclusions give the full strip and permutation inputs.

\cref{r10.9,r10.10} apply. The obstruction in
\cref{r12.3}, the balanced covers of \cref{r2.4}, and \cref{r12.6} give the
saturation in \cref{r12.8}. The six-cover structure from
\cref{r13.3} is excluded by \cref{r14.5}, proving the alternative.
\end{proof}

\subsection{Proof of the upper bound}\label{sec:15}
\begin{proof}[Proof of \cref{r1.1}]
Fix the transfer tolerance, extraction hierarchy and regularity
parameters, and take $n$ above all resulting finite thresholds.
A three-colouring of $K_{3n,3n}$ without a monochromatic $C_{2n}$
gives the admissible profile in \cref{r2.2}. Its three possible
configurations in \cref{struct:alternative} are excluded by
\cref{finite:strip,r9.7,r10.8}, respectively. This proves the
upper bound. Together with \cref{r1.2}, it yields
$\br(C_{2n})=3n$ for all sufficiently large $n$.
\end{proof}

\section*{Declaration on the use of generative AI}
The lower-bound construction is entirely due to the authors.
For the upper bound, the authors supplied the proof framework and directed
its revision. OpenAI's GPT-6.1 Sol and Anthropic's Claude (Fable 5.1) assisted with checking intermediate arguments, completing local
derivations, and refining the multi-core capacity argument and the
thresholds before and after cleaning. They also assisted with English
editing. The authors are responsible for the mathematical claims and the
final text.

\end{document}